\documentclass[12pt,a4paper]{article}

\usepackage[a4paper,lmargin={2.5cm},rmargin={2.5cm},tmargin={2.5cm},bmargin={2.5cm}]{geometry}

\usepackage{amsmath}
\usepackage{amssymb}
\usepackage{amstext}
\usepackage{amsthm}

\newtheorem{definition}{Definition}
\newtheorem{notation}{Notation}

\newcommand{\dd}{\mathrm d}
\newcommand{\refel}{\mathrm{ref}}
\newcommand{\recon}{\mathrm{recon}}
\newcommand{\cfl}{\mathrm{cfl}}

\usepackage{graphicx}
\usepackage[dvipsnames]{xcolor}
\usepackage{scalerel}

\usepackage{tikz}

\usetikzlibrary{shapes.geometric}
\usetikzlibrary{shapes.misc}

\newcommand{\symtri}{\protect\tikz \protect\node[isosceles triangle, draw=black, anchor=south, rotate=135, isosceles triangle apex angle=90, minimum height=0.12cm, inner sep=0, outer sep=0] at (0, 0) {};} 

\usepackage{array}
\usepackage{hhline}

\usepackage[english]{babel}
\usepackage{authblk}

\usepackage{csquotes}

\usepackage{hyperref}
\usepackage{subcaption}

\title{A generalized Active Flux method of arbitrarily high order in two dimensions}
\author[1]{Wasilij Barsukow}
\author[2]{Praveen Chandrashekar}
\author[3]{Christian Klingenberg}
\author[3]{Lisa Lechner}
\affil[1]{Institut de Mathématiques de Bordeaux (IMB), CNRS UMR 5251, 351 Cours de la Libération, 33405 Talence, France}
\affil[2]{Centre for Applicable Mathematics, Tata Institute of Fundamental Research, Bengaluru, 560065, India}
\affil[3]{University of Würzburg, Institute of Mathematics, Emil-Fischer-Straße 40, 97074 Würzburg, Germany}
\date{\vspace{-15mm}}

\begin{document}

\maketitle

\begin{abstract}
The Active Flux method can be seen as an extended finite volume method. The degrees of freedom of this method are cell averages, as in finite volume methods, and in addition shared point values at the cell interfaces, giving rise to a globally continuous reconstruction. Its classical version was introduced as a one-stage fully discrete, third-order method. Recently, a semi-discrete version of the Active Flux method was presented with various extensions to arbitrarily high order in \emph{one} space dimension. 
In this paper we extend the semi-discrete Active Flux method on \emph{two}-dimensional Cartesian grids to arbitrarily high order, by including moments as additional degrees of freedom (hybrid finite element--finite volume method). The stability of this method is studied for linear advection. For a fully discrete version, using an explicit Runge-Kutta method, a CFL restriction is derived. We end by presenting numerical examples for hyperbolic conservation laws.


\emph{Mathematics Subject Classification (MSC2020):}
65M08, 
65M12, 
65M20, 
35L65.

\end{abstract}

\section{Introduction}

The solution of systems of hyperbolic conservation laws is required by many different physical problems from fluid dynamics. 
As most of these problems are fairly complex, numerical methods are needed to solve them. Finite volume methods evolve cell averages and maintain conservation in a discrete sense which is important for convergence to a weak solution. Godunov's method \cite{God1959}, for example, is a seminal approach to hyperbolic conservation laws. More recently, Roe and coauthors introduced the Active Flux method \cite{ER2011a,ER2011b,ER2013} inspired by van Leer's scheme V \cite{VanLeer1977} for linear advection. The new method not only uses cell averages but also point values at the cell interfaces. This allows to reconstruct the solution piecewise parabolically in a globally continuous fashion. In \cite{VanLeer1977}, the point values at the cell interfaces were updated by tracing back the characteristics. This way, a method of third order for linear advection in one spatial dimension arose. Roe's one-stage fully discrete Active Flux method (which we will call the \emph{classical Active Flux method}, for a review on Cartesian grids see e.g. \cite{CH2023}) extends this idea to certain nonlinear and/or multi-dimensional hyperbolic problems by replacing characteristic transport with suitable either exact or approximate evolution operators (e.g. \cite{maeng17,fan17,BHKR2019,Bar2020,chudzik24}). Depending on the system of conservation laws finding such an evolution operator may be quite challenging. Recently the authors in \cite{Abgrall2023, AB2023ExtensionAF, AB2023FEFV, ABK2025} proposed a semi-discrete formulation of the Active Flux method (which we will call the \emph{generalized Active Flux method}) that allows for a more flexible approach to update the point values. 

In \cite{AB2023ExtensionAF}, several extensions of this method to arbitrary order in one spatial dimension are presented, while the multi-dimensional method in \cite{ABK2025} was only third-order accurate. In this paper we extend this work to obtain a semi-discrete generalized Active Flux method of arbitrarily high order on Cartesian grids in two spatial dimensions. The design of the method is based on a hybrid finite element--finite volume method (see also \cite{AB2023FEFV}) and also takes inspiration from the one-dimensional arbitrary-order Active Flux method with additional point values developed in \cite{AB2023ExtensionAF}.

While this paper focuses on Cartesian grids, we want to also mention \cite{ALL2025}, where a semi-discrete Active Flux method on triangular meshes is presented and results for a third- and fourth-order accurate method are shown, using the cell average and point values at the cell interfaces. 
Furthermore, \cite{He2021} and \cite{Roe2021,RS2023,Samani2024} consider extensions of the (classical) fully discrete Active Flux method to arbitrary order. The first approach uses additional point values. In the latter, they typically enrich the stencil using additional information in the form of derivatives at the already defined point values at the cell interfaces (also called \enquote{Hermite Active Flux}).

Here, we aim at using higher moments in the cells and additional point values at the cell interfaces. For efficiency, we especially aim at achieving arbitrarily high order with a minimal number of degrees of freedom while maintaining the compact stencil. Another possibility is a tensor-like extension of the hybrid finite element--finite volume method in \cite{AB2023FEFV}, which is also formally defined in \cite{MthJK2024}, but it uses more degrees of freedom than actually necessary for the required order of accuracy.

Using our approach, key features of Active Flux like global continuity of the reconstructed solution and a compact stencil of the method in space are retained. As the method is based on the semi-discrete approach, it has the potential to be more easily applicable to different systems of conservation laws than the classical Active Flux. For integration in time we use an explicit Runge-Kutta method.
While classical Active Flux is a one-stage fully discrete method, the multiple stages of a Runge-Kutta method widen the effective stencil of the generalized Active Flux method. In \cite{Roe2021}, this has been linked to a reduction of the maximal CFL number, compared to the classical Active Flux method. Yet, \cite{BKKL2025} shows encouraging similarities of the classical and generalized Active Flux method of third order on Cartesian grids for linear acoustics, such as stationarity preservation (see \cite{BHKR2019} for stationarity preservation of the classical Active Flux).

The paper is structured as follows: In Section \ref{Chap_ReviewSemiDiscreteAF} we briefly recall arbitrarily high-order one-dimensional Active Flux from \cite{AB2023ExtensionAF,AB2023FEFV} and third-order two-dimensional Active Flux from \cite{ABK2025}. Then, our high order method is defined in Section \ref{Chap_GenAF2DCartesian}. Section \ref{Chap_StabilityAnalysis} focuses on a stability analysis for linear advection of both the semi-discrete and the fully discrete method using RK3. Numerical results up to seventh order are presented for our high-order method in Section \ref{Chap_NumEx}, including an example for the 2-d compressible Euler equations.

\section{An overview of the semi-discrete Active Flux method}
\label{Chap_ReviewSemiDiscreteAF}
Consider the hyperbolic system of conservation laws
\begin{equation}\label{Eq_SysCL}
\partial_t q(t, \mathbf x) + \nabla  \cdot \mathbf f(q(t,\mathbf x)) = 0
\end{equation}
for \(q\colon \mathbb R_0^+ \times \Omega \rightarrow \mathbb R^s\) on the domain \(\Omega \subseteq \mathbb R^d\). 
Throughout this paper the focus is on the two-dimensional case \(d=2\) with \(\mathbf f= (f^x, f^y)\), \(f^x, f^y\colon\mathbb R^s \rightarrow \mathbb R^s\). The notation used here is similar to \cite{ABK2025}.

Furthermore, we restrict ourselves to rectangular domains \(\Omega = [x_\mathrm{min}, x_\mathrm{max}] \times [y_\mathrm{min}, y_\mathrm{max}]\) discretized with Cartesian grids with equidistant cell widths $\Delta x$ and $\Delta y$ in $x$- and $y$-direction. For a Cartesian grid of size $N_x \times N_y$ the grid cells will be denoted as \(C_{ij} = [x_{i-\frac{1}{2}}, x_{i+\frac{1}{2}}] \times [y_{j-\frac{1}{2}}, y_{j+\frac{1}{2}}]\) for \(i = 0, \dots, N_x-1, j = 0, \dots N_y-1\) centered at \( \mathbf x_{ij} = (x_i, y_j)\).\\

Next, we will give an overview of semi-discrete Active Flux methods, that we will use and extend in this paper.

\subsection{A third-order method on two-dimensional Cartesian grids}
\label{Chap_GenAF3rd2dCartesian}
\begin{figure}
\centering
	\begin{subfigure}[b]{0.45\textwidth}
	\centering		
		\includegraphics[scale=1.4]{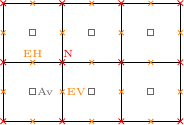}
		\caption{Sketch of degrees of freedom: cell averages (Av) and shared point values at vertices (N) and horizontal (EH) and vertical (EV) edges.}
		\label{Fig_SketchGrid2D3rd} 
	\end{subfigure}
	\quad
	\begin{subfigure}[b]{0.45\textwidth}
	\centering
		\includegraphics[scale=0.5]{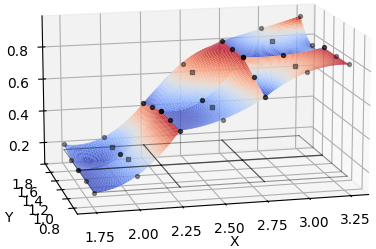}
		\caption{Visualization of the concept for the globally continuous reconstruction via a constructed example.}
		\label{Fig_ExampleReconstruction2D3rd} 
	\end{subfigure}
	\caption{Degrees of freedom and reconstruction for a third-order Active Flux method on two-dimensional Cartesian grids.}
\end{figure}

In \cite{ABK2025} a semi-discrete Active Flux method of third order on Cartesian grids was presented for the two-dimensional case. 
As for the classical Active Flux method on Cartesian grids (see e.g. \cite{BHKR2019,Helzel2019}),
the degrees of freedom for each cell \(C_{ij}\) are the cell average
\begin{equation}\label{Eq_CellAvAF3rd}
	\bar q_{ij}(t) = \frac{1}{\Delta x \Delta y} \int_{C_{ij}} q(t,\mathbf x) \mathrm d \mathbf x
\end{equation}
and the nodal, vertical and horizontal point values
\begin{equation}\label{Eq_NodesAF3rd}
	q_{i+\frac{1}{2}, j+ \frac{1}{2}}(t) = q(t, x_{i+\frac{1}{2}}, y_{j+ \frac{1}{2}}), 
	\quad q_{i, j+ \frac{1}{2}}(t) = q(t, x_{i}, y_{j+ \frac{1}{2}}), 
	\quad q_{i+\frac{1}{2}, j}(t) = q(t, x_{i+\frac{1}{2}}, y_{j})
\end{equation}
at the cell interfaces. All point values are shared, i.e. each cell has access to eight point values at $\mathbf x_p$, $p \in \{(i \pm \tfrac12, j\pm \tfrac12), (i \pm \tfrac12,j), (i,j \pm \tfrac12)\}$ (see Figure \ref{Fig_SketchGrid2D3rd}).

The reconstruction on cell \(C_{ij}\) is the biparabolic polynomial \(q_{ij, \mathrm{recon}}: [-\tfrac{\Delta x}{2}, \tfrac{\Delta x}{2}]\times[-\tfrac{\Delta y}{2},\tfrac{\Delta y}{2}] \rightarrow \mathbb R^s\) (see \cite{ABK2025, BHKR2019}) 
defined on a reference cell with \(\mathbf x_\mathrm{ref} := \mathbf x - \mathbf x_{ij}\) and \(\mathbf x_{\mathrm{ref},p} := \mathbf x_p - \mathbf x_{ij}\) which has to satisfy
\begin{align*}
	&\frac{1}{\Delta x \Delta y} \int_{-\frac{\Delta y}{2}}^{\frac{\Delta y}{2}}\int_{-\frac{\Delta x}{2}}^{\frac{\Delta x}{2}} q_{ij, \mathrm{recon}}(\mathbf x_\mathrm{ref}) \mathrm d x_\mathrm{ref} \mathrm d y_\mathrm{ref} = \bar q_{ij},\\
	&q_{ij, \mathrm{recon}}(\mathbf x_{\mathrm{ref},p}) = q_p \quad \forall p \in \{(i \pm \tfrac12, j\pm \tfrac12), (i \pm \tfrac12,j), (i,j \pm \tfrac12)\}.
\end{align*}
It can also be written in terms of shape functions, see e.g. \cite{BKKL2025}. This gives rise to a globally continuous reconstruction \(q_\mathrm{recon}\) on \(\Omega\) (see Figure \ref{Fig_ExampleReconstruction2D3rd}) with 
\begin{equation}
	q_\recon|_{C_{ij}}(\mathbf x) := q_{ij, \recon}(\mathbf x_\refel).
\end{equation}

For the update of the point values \(q_p, p \in \{(i+\tfrac{1}{2},j+\tfrac{1}{2}), (i+\tfrac{1}{2},j), (i,j+\tfrac{1}{2})\}\) the quasi-linear form of \eqref{Eq_SysCL} at $\mathbf x_p$
\begin{equation*}
	\frac{\mathrm d}{\mathrm d t}q_p(t)  +  Df^x(q_p(t)) \partial_x q\vert_p(t) +  Df^y(q_p(t)) \partial_y q\vert_p(t)  = 0 
\end{equation*}
is considered.
In order to introduce stabilization via upwinding, the Jacobians \(J^x = Df^x(q_p), J^y = Df^y(q_p) \) are split according to the sign of the wave speeds
\begin{align*}
	(J^x)^{\pm} &:= T^x\mathrm{diag}(\lambda_{x,1}^{\pm},\dots , \lambda_{x,s}^{\pm})(T^x)^{-1},\\
	\qquad \lambda_{x,k}^{+} &:= \max(0, \lambda_{x,k}), \quad \lambda_{x,k}^{-} := \min(0, \lambda_{x,k}) \qquad \forall k = 1, \dots, s. 
\end{align*}
The derivatives \(\partial_x q\vert_p\) are approximated with finite difference formulas in upwind direction \((D_x^{\pm})_p q\) and \((D_y^{\pm})_p q\) derived from the reconstruction \(q_\mathrm{recon}\). 
In directions tangential to the cell interface for the horizontal and vertical edge point, the derivative of the reconstruction is continuous and no upwinding is included. \\

The update of the cell averages is obtained by integrating the conservation law \eqref{Eq_SysCL} over the cell \(C_{ij}\) and applying the divergence theorem
\begin{equation*}
	\frac{\mathrm d}{\mathrm d t} \bar q_{ij} (t) + \frac{1}{\Delta x \Delta y} \int_{\partial C_{ij}} \mathbf f\left(q(t, \mathbf x)\right) \cdot \mathbf n \mathrm{d}S = 0.
\end{equation*}
The reconstructed solution is continuous at the cell interfaces and the flux terms over the cell interfaces are numerically integrated with a Gauss-Lobatto quadrature with three points, which coincide with the point values at the cell interfaces.

\subsection{Methods of arbitrary order in one dimension}
In \cite{AB2023ExtensionAF}, different methods to extend semi-discrete Active Flux to arbitrarily high order in one dimension were suggested. Here, two of these concepts are reviewed, and a variation of the second method is presented.

\subsubsection{A method with higher moments}
\label{Chap_GenAF1dMoments}
In this paper, we will mainly focus on the hybrid finite element--finite volume method as also defined in \cite[Definition 4.1 (Method A)]{AB2023FEFV}. 
For a method of order $N+1$ on a compact stencil, $N+1$ degrees of freedom accessible to a cell are needed to approximate a polynomial of degree $N$. 
In addition to the traditional degrees of freedom, i.e. the cell averages and the point values at the cell interfaces, this method uses higher moments 
\begin{equation*}
	q_i^{(k)}(t) := A_k \int_{x_{i-\frac{1}{2}}}^{x_{i+\frac{1}{2}}} q(t, x) b_k(x-x_i) \mathrm d x  \qquad k \in \mathbb N.
\end{equation*}
The test functions for the moments are \(b_k: [-\frac{\Delta x}{2},\frac{\Delta x}{2}]\rightarrow \mathbb R, b_k(x) = x^k\) with the normalization factors \(A_k = \tfrac{(k+1)2^k}{\Delta x^{k+1}}\). This allows to obtain more information on a compact stencil. 
The update of the cell average $\bar q_i = q_i^{(0)}$ and higher moments is given as
\begin{align*}
	\frac{\mathrm d}{\mathrm d t} q_i^{(k)}(t) = &-A_k\left(f(q_{i+\frac{1}{2}}(t)) b_k\left(\tfrac{\Delta x}{2}\right) 
											- f(q_{i-\frac{1}{2}}(t))b_k\left(-\tfrac{\Delta x}{2}\right)\right)\\
								         &+ A_k \int_{x_{i-\frac{1}{2}}}^{x_{i+\frac{1}{2}}} f(q(t, x)) b'_k(x-x_i) \mathrm d x \quad \forall k = 0,\dots, N-2
\end{align*}
for a method of order \(N+1\). The update of the point values \(q_{i+\frac{1}{2}} = q(x_{i+\frac{1}{2}})\) uses a Jacobian splitting of $Df(q)$ to include upwinding
\begin{equation*}
	\frac{\mathrm d}{\mathrm d t} q_{i+\frac{1}{2}}(t) = -Df(q_{i+\frac{1}{2}})^+ (D_x^+)_{i+\frac{1}{2}}q - Df(q_{i+\frac{1}{2}})^- (D_x^-)_{i+\frac{1}{2}}q.
\end{equation*}
The polynomial reconstructions \(q_{i, \mathrm{recon}}:[-\frac{\Delta x}{2}, \frac{\Delta x}{2}] \rightarrow \mathbb R^s\) of degree $N \geq 2$ are used to find a finite difference approximation for \((D_x^\pm)_{i+\frac{1}{2}}q\).

\subsubsection{A method with additional point values}
\label{Chap_AFaddPV1d}
Another possibility to achieve arbitrarily high order of the Active Flux method is to introduce further point values \(q_{i,\xi_a}= q(x_i + \Delta x \xi_a), \xi_a \in (-\tfrac{1}{2}, \tfrac{1}{2})\) within each cell in addition to the cell averages and the point values at the cell interfaces \(q_{i\pm\frac{1}{2}}\). This method will especially be relevant for the discussion of stability in Section \ref{Chap_StabilityAnalysis}.
To achieve a method of order \(N+1\), \(N-2\) additional point values are needed and similarly to Section \ref{Chap_GenAF1dMoments}, a polynomial reconstruction can be defined on the cells \(C_i\) (except in the case of a symmetric distribution of an odd number of point values). 
In \cite{AB2023ExtensionAF} this method has been introduced using evolution operators for the point value updates.\\

Before proceeding with the arbitrarily high-order two-dimensional method, we propose a variation of the above method here. Since it is a challenge to find exact or high-order approximate evolution updates for the point values for general systems of conservation laws and having in mind the point value updates of the method reviewed in Section \ref{Chap_GenAF1dMoments}, it seems straightforward to write a variation of the method from \cite{AB2023ExtensionAF}, for $N+1\geq 4$, as the following semi-discretization:
\begin{align*}
	\left\lbrace
	\begin{aligned}
		\frac{\mathrm d}{\mathrm d t} \bar q_i(t) &= -\frac{f(q_{i+\frac{1}{2}}(t)) - f(q_{i-\frac{1}{2}}(t))}{\Delta x}\\
		\frac{\mathrm d}{\mathrm d t} q_{i+\frac{1}{2}}(t) &= -F\Big(q_{i-\frac{1}{2}}(t), q_{i,\xi_0}(t), \dots, q_{i,\xi_{N-3}}(t),
												q_{i+\frac{1}{2}}(t), \\
												& \hspace{6em} q_{i+1,\xi_0}(t), \dots, q_{i+1,\xi_{N-3}}(t),
												 q_{i+\frac{3}{2}}(t)\Big)\\
		\frac{\mathrm d}{\mathrm d t} q_{i,\xi_a}(t) &= -\tilde F_a\Big(q_{i-\frac{1}{2}}(t), q_{i,\xi_0}(t), \dots, q_{i,\xi_{N-3}}(t),
												q_{i+\frac{1}{2}}(t)\Big) \quad \forall a \in \{0, \dots, N-3\}
	\end{aligned}
	\right.
\end{align*}
with consistent approximations $F$, $\tilde F_a$ of $\partial_x f(q)$ at $x_{i+\frac{1}{2}}$ and $x_i + \Delta x \xi_a$, respectively.
As in Section \ref{Chap_GenAF1dMoments}, a Jacobian splitting can be considered and a finite difference formula can be derived from the reconstructed polynomial. For the point values at the cell interfaces it is natural to use the reconstruction of the cell in upwind direction such that \((D^+_x)_{i+\frac12}q := \frac{\dd}{\dd x_\mathrm{ref}}q_{i, \mathrm{recon}}(x_\mathrm{ref})\vert_{x_\mathrm{ref}=\frac{\Delta x}{2}}\),  \((D^-_x)_{i+\frac12}q := \frac{\dd}{\dd x_\mathrm{ref}}q_{i+1, \mathrm{recon}}(x_\mathrm{ref})\vert_{x_\mathrm{ref}=-\frac{\Delta x}{2}}\). For the point values within the cell, the reconstruction is uniquely defined and \((D_x)_{i, \xi_a}q := \frac{\dd}{\dd x_\mathrm{ref}}q_{i, \mathrm{recon}}(x_\mathrm{ref})\vert_{x_\mathrm{ref}=\Delta x \xi_a}\). This resembles the approximation of the derivative tangential to a cell interface in Section \ref{Chap_GenAF3rd2dCartesian}. The version described in this Section will serve as inspiration in Section \ref{Chap_StabilityAnalysis}.

\section{A semi-discrete Active Flux method of arbitrarily high order on two-dimensional Cartesian grids}
\label{Chap_GenAF2DCartesian}
In the following Section we discuss the extension of the semi-discrete, third-order Active Flux method (as recapitulated in Section \ref{Chap_GenAF3rd2dCartesian}) to arbitrary order of accuracy $N+1$. It will be constructed in such a way that 
\begin{itemize}
	\item the reconstructed numerical solution \(q_\mathrm{recon}\) is globally continuous,
	\item it uses a compact stencil with a minimal number of degrees of freedom per cell, and
	\item the cell average is always included among the degrees of freedom.
\end{itemize}
As for the third-order method, the shared point values at the cell interfaces will guarantee that the reconstruction is globally continuous and the inclusion of the cell averages will ensure that the resulting method is conservative. To gain sufficient information for the reconstruction on a compact stencil, higher moments are defined. Hence, this work extends the hybrid finite element--finite volume approach (see Section \ref{Chap_GenAF1dMoments}) from one spatial dimension to two spatial dimensions on Cartesian grids.\\
The way the point values are distributed on the edge has an influence on stability as shall be discussed in Section \ref{Chap_StabAnaSemiDiscreteGenAF}. There, similar considerations are made as for the Active Flux method with additional point values (see Section \ref{Chap_AFaddPV1d}).

\noindent\emph{Remark:} This approach to extend the generalized Active Flux method is only one possibility. 
It uses point values at the cell interfaces analogously to the third order method in 2-d and moments in the interior of the cell, extending the 1-d method from Section \ref{Chap_GenAF1dMoments}. Using moments in the interior gives a natural way to minimize the number of degrees of freedom. Another possibility can be a tensor-like extension of the 1-d method with moments which will also briefly be considered here. 
Further extensions can be thought of. For example the method described in Section \ref{Chap_AFaddPV1d} could be extended. A similar method considering point values only is the face-upwinded spectral element (FUSE) method \cite{PP2024} on quadrilateral meshes.

\subsection{Reference element for Cartesian grids}
From now on, we define the reference element as \(C^\mathrm{ref} := \left[-\frac12, \frac12\right]\times\left[-\frac12,\frac12\right]\) and the transformation \(\Phi_{ij}: C_{ij} \rightarrow C^\mathrm{ref}\)
\begin{equation}
\Phi_{ij}(x,y) = \left(\begin{array}{c} \frac{x-x_i}{\Delta x}\\
						      \frac{y-y_j}{\Delta y}\end{array}\right)
\end{equation}
with its inverse \(\Phi_{ij}^{-1}: C^\mathrm{ref} \rightarrow C_{ij}\) 
\begin{equation}
\Phi_{ij}^{-1}(x_\mathrm{ref},y_\mathrm{ref}) = \left(\begin{array}{c} {x_i + \Delta x x_\mathrm{ref}}\\
												{y_j + \Delta y y_\mathrm{ref}}\end{array}\right).
\end{equation}

\subsection{Degrees of freedom}
\label{Chap_DOFs}
\begin{figure}
\centering
	\includegraphics[scale=1.3]{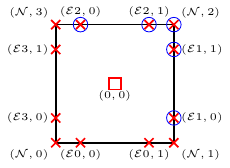}
	\caption{Superscripts specifying the degrees of freedom for a cell $C_{ij}$ for GenAF$(M_{\symtri})$ of fourth order (a generalized Active Flux method with minimal number of degrees of freedom, see Notation \ref{Not_GenAFMMin}). The circled point values mark the points \emph{belonging} to a cell.}
	\label{Fig_GenAFDOFsIndex} 
\end{figure}

The degrees of freedom used for the generalized Active Flux method of order $N+1$ are shared point values at the cell interfaces and moments.
The notation we introduce here uses a subscripted index $ij$ to identify the cell and a superscripted index to specify the degree of freedom. An example for the superscripts can be seen in Figure \ref{Fig_GenAFDOFsIndex}.\\

As described above, the point values at the cell interfaces are chosen such that the reconstruction \(q_\mathrm{recon}\) is globally continuous on the domain \(\Omega\) and maximal use of these shared degrees of freedom is made. 
\(N+1\) shared point values at each cell interface are required to define a unique polynomial of degree $N$ along each cell interface and achieve \(q_\mathrm{recon} \in C^0(\Omega, \mathbb R^s)\). 
First, the four nodes (= vertices) are included. This leaves $N-1$ points per edge, i.e. \(4(N-1)\) edge points in total. 
For cell \(C_{ij}\) the point values are defined as 
\begin{equation}\label{Eq_GenAF2DarbPV}
	q_{ij}^{(p)}(t):= q\big(t, \mathbf x_{ij}^{(p)}\big) \quad \forall p \in I_{\mathcal P}
\end{equation}
with \(\mathbf x_{ij}^{(p)} \in \partial C_{ij}\). 
The index set 
\begin{equation}\label{Eq_IndexPointsAccCell}
	I_{\mathcal P} := I_{\mathcal N} \cup I_{\mathcal E}
\end{equation}
defines all points \emph{accessible} to the cell. These are 
\begin{itemize}
\item the nodes \(\big\{\mathbf x_{ij}^{(n)}\big\}_{n \in I_{\mathcal N}}\), \(I_{\mathcal N} := \{(\mathcal N,0),\dots, (\mathcal N,3)\}\) located at
\begin{equation*}
	(x_{i \pm \frac{1}{2}}, y_{j\pm \frac{1}{2}})
\end{equation*}
and assembled counterclockwise starting with the lower left node and 
\item the edge points \(\big\{\mathbf x_{ij}^{(e, a)}\big\}_{(e, a) \in I_{\mathcal E}}\), \(I_{\mathcal E} := \{(e,a) | e = \mathcal E0,\dots, \mathcal E3, a = 0,\dots, N-2\}\) with the edges assembled counterclockwise. The horizontal edge points with \(e = 0, 2\) are located at 
\begin{equation*}
	(x_{i +\xi_{a}},y_{j \pm \frac{1}{2}})=(x_i +\xi_{a} \Delta x,y_{j \pm \frac{1}{2}})
\end{equation*}
and the vertical edge points with \(e = 1, 3\) at 
\begin{equation*}
	(x_{i \pm \frac{1}{2}},y_{j +\xi_{a}}) = (x_{i \pm \frac{1}{2}},y_j  +\xi_{a} \Delta y)
\end{equation*}
with \(\xi_{a} \in (-\tfrac{1}{2}, \tfrac{1}{2})\). The choice of \(\xi_{a}\) is discussed in Section \ref{Chap_StabAnaSemiDiscreteGenAF}. 
\end{itemize}

Since the point values are shared between adjacent cells it is possible to define the points \emph{belonging} to a cell, e.g. chosen as the upper right node \(n = 2\) as well as the edge points on the upper horizontal \(e = \mathcal E2\) and the right vertical edge \(e = \mathcal E1\), assembled in the index set 
\begin{equation}\label{Eq_IndexPointsBelCell}
	\widetilde {I_{\mathcal P}} := \{(\mathcal N, 2)\} \cup \{(e,a) | e= \mathcal E1,\mathcal E2, a = 0,\dots, N-2\}.
\end{equation}
Thus, on average only $1$ node and $2(N-1)$ edge points have to be updated per cell. \\

The moments of order \(m = k + l\), $k,l \geq 0$ for cell \(C_{ij}\) are defined as
\begin{equation}\label{Eq_GenAF2DarbMoments}
q_{ij}^{(k,l)}(t) := A_{k,l} \int_{C_{ij}} b_{k,l}(\Phi_{ij}(\mathbf x)) q(t, \mathbf x) \mathrm d \mathbf x \quad \forall (k,l) \in I_{\mathcal M}
\end{equation}
with the test functions \((b_{k,l})_{k,l \in \mathbb N_0}\)
\begin{equation}\label{Eq_GenAFTestFunctionMoments}
b_{k,l}:  \left[-\tfrac{1}{2}, \tfrac{1}{2}\right] \times \left[-\tfrac{1}{2}, \tfrac{1}{2}\right] \rightarrow \mathbb{R}, \quad b_{k,l}(x,y) = x^ky^l
\end{equation}
and a corresponding normalization factor
\begin{equation}\label{Eq_GenAFNormalizationFactorMoments}
A_{k,l} = \frac{(k+1)2^k(l+1)2^l}{\Delta x \Delta y}.
\end{equation}
Thus, there are \(m+1\) moments of order \(m\) and the cell averages \(\bar q_{ij}\) are equal to the moments \(q_{ij}^{(0,0)}\), i.e. $m=0$.\\
Beyond the requirement of including the average among the degrees of freedom, there are different possibilities to choose the moments. The two choices considered here are 
\begin{align}
I_{\mathcal M} &= I_{\mathcal M,\square} := \{(k,l) | k,l \in \mathbb N_0, 0 \leq k,l \leq N-2\} && \text{(tensor-like)} \label{eq:momentstensorlike}\\
I_{\mathcal M} &= I_{\mathcal M,\symtri} := \{(k,l) | k,l \in \mathbb N_0, 0 \leq k+l \leq \max\{0, N-4\}\} && \text{(triangle-like)} \label{eq:trianglelike}
\end{align}
for $N\geq2$. They entail different reconstruction spaces which will be discussed in the following Section.\\

\begin{figure}
\centering
	\begin{subfigure}[b]{0.25\textwidth}
	\centering		
		\includegraphics[scale=1.1]{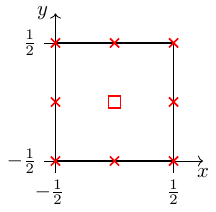}
		\caption{\(N+1 = 3\).}
		\label{Fig_RefEl3rd} 
	\end{subfigure}
	\begin{subfigure}[b]{0.25\textwidth}
	\centering
		\includegraphics[scale=1.1]{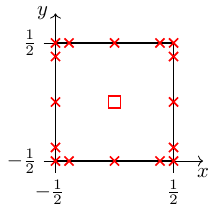}
		\caption{\(N+1 = 5\).}
		\label{Fig_RefEl5th} 
	\end{subfigure}
	\begin{subfigure}[b]{0.25\textwidth}
	\centering
		\includegraphics[scale=1.1]{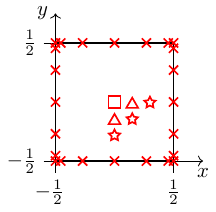}
		\caption{\(N+1 = 7\).}
		\label{Fig_RefEl7th} 
	\end{subfigure}
	\begin{subfigure}[b]{0.2\textwidth}
	\centering
		\includegraphics[scale=1, trim=-0.5cm -1.5cm 0 0]{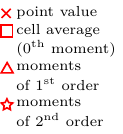}
	\end{subfigure}
	\caption{Reference elements \(C^\mathrm{ref}\) for GenAF$(M_{\symtri})$ of order $N+1$ (a generalized Active Flux method with minimal number of degrees of freedom, see Notation \ref{Not_GenAFMMin}).}
	\label{Fig_GenAFMinDOFsRefEl} 
\end{figure}

The point values and moments define all our degrees of freedom and the index set of the nodal and modal degrees of freedom accessible to each cell is denoted by
\begin{equation}
	I_\mathrm{dof} := I_{\mathcal P} \cup I_{\mathcal M}.
\end{equation}
and analogously we denote by
\begin{equation}
	\widetilde{I_\mathrm{dof}} := \widetilde{I_{\mathcal P}} \cup I_{\mathcal M}
\end{equation}
the index set of all degrees of freedom belonging to a cell. 
Figure \ref{Fig_GenAFMinDOFsRefEl} shows some examples for the reference elements of spatial order \(N+1\) using the triangle-like moments  \eqref{eq:trianglelike}. 
The updates of the degrees of freedom are discussed in the Sections \ref{Chap_GenAF2DarbPV} and \ref{Chap_GenAF2DarbMoments}.\\

\noindent\emph{Remark:} Finite elements using nodal and modal degrees of freedom are also called hybrid finite elements (see e.g. \cite[Chapter~6.3.3,Chapter~7.6]{ErnGuermond2021}).

\subsection{Reconstruction}
\label{Chap_GenAF2DReconstruction}
Next, the spatial reconstruction \(q_\mathrm{recon}: \Omega_h \rightarrow \mathbb R^s\), \(\Omega_h:= \cup_{i,j} C_{ij}\), is considered (for 1-d it reduces to the reconstruction in \cite{AB2023FEFV}, and for third order in 2-d to the one in \cite{BHKR2019}).
The polynomial solution reconstruction on cell \(C_{ij}\) 
\begin{align}\label{Eq_DefRecon1}
	q_\mathrm{recon}|_{C_{ij}}&\colon C_{ij} \rightarrow \mathbb R^s & q_\mathrm{recon}|_{C_{ij}} & \in (P^\mathrm{recon})^s
\end{align}
has to fulfill \eqref{Eq_GenAF2DarbPV} and \eqref{Eq_GenAF2DarbMoments}, 
and for the method to be of order $N+1$ in space it must at least lie in the bivariate polynomial space 
\begin{align}
P^N:=\mathrm{span}\{x^m y^n | m, n \in \mathbb N_0, 0 \leq m+n \leq N\}
\end{align}
of maximal total degree $N$, for each element, i.e. \(P^N \subset P^\mathrm{recon}\). 
Here, we use all degrees of freedom \(q_{ij}^{(r)}, r \in I_\mathrm{dof}\) accessible to the cell to define \eqref{Eq_DefRecon1}. Clearly, a necessary condition for unisolvence is that the dimension $|P^\mathrm{recon}|$ of $P^\mathrm{recon}$ matches the number of degrees of freedom. 
Below, two possible choices for the reconstruction space are presented. Later on, they are used to define a finite element for the generalized Active Flux method.

\subsubsection{Tensor-like reconstruction space}
As a natural extension from the one-dimensional case (see \cite{AB2023FEFV,AB2023ExtensionAF}) the bi-polynomial space
\begin{equation}
	P^{N,N}:= \mathrm{span}\{x^m y^n | m,n \in \{0, \cdots, N\}\}
\end{equation}
can be considered for \(P^\mathrm{recon}\), similar to tensor-product finite elements (see e.g. \cite[Chapter~6.4]{ErnGuermond2021}).

This choice satisfies \(P^N \subset P^\mathrm{recon}\) and requires \(|P^\mathrm{recon}|=(N+1)^2\) degrees of freedom. As discussed in Section \ref{Chap_DOFs}, $4N$ of these must be point values which leaves $(N-1)^2$ modal degrees of freedom within each cell, i.e. the cell average and the higher moments. Thus, the tensor-like reconstruction space is associated to the tensor-like choice \eqref{eq:momentstensorlike} of moments spanning \(P^{N-2,N-2}\). This choice reduces to the reconstruction space for the third-order Active Flux method as first introduced in \cite{BHKR2019,Helzel2019}, see also Section \ref{Chap_GenAF3rd2dCartesian}.\\

For the resulting method we introduce the following notation:
\begin{notation}
\label{Not_GenAFMTens}
	The \emph{tensor-like} generalized Active Flux method on two-dimensional Cartesian grids with $P^\mathrm{recon}=P^{N,N}$ shall be denoted by \emph{GenAF$(M_\square)$}.
\end{notation}

\subsubsection{Reconstruction space with minimal number of degrees of freedom}
High-order methods need to evolve many degrees of freedom. We thus aim at finding a computationally more efficient choice for \(P^\mathrm{recon}\) with only the minimal possible number of degrees of freedom needed to achieve order \(N+1\) on a Cartesian cell.

Setting \(P^\mathrm{recon} = P^N\) turns out not to be sufficient for unisolvence, as discussed below. Instead we define
\begin{align}\label{Eq_PReconMin}
	P^\mathrm{recon} :=
	P^\mathrm{recon}_{\min} := 
	\left\lbrace	
	\begin{aligned}
	&S^N \oplus \mathrm{span}\{x^2y^2\} &&\text{ for } N = 2,3\\	
	&S^N &&\text{ for } N \geq 4.
	\end{aligned}
	\right.
\end{align}
using the so called serendipity space (see also Table \ref{Table_SerendipitySpaceBasis})
\begin{equation}
S^N:=P^N \oplus \mathrm{span}\{x^N y, x y^N\}
\end{equation}
known from serendipity finite elements (see e.g. \cite[Chapter~6.4.3]{ErnGuermond2021} and \cite{AA2011}), which are based on the idea of reducing the degrees of freedom from the tensor-product finite element.
\begin{table}
\fontsize{10pt}{10pt}\selectfont
\newcommand{\so}{\scaleobj{0.6}}
\newcommand{\rb}{\reflectbox}
\centering
	\begin{tabular}{| c | c !{\vrule width 1pt} c | c | c | c | c |}
	\hline
	\(1\) 	\phantom{\Big|}		& \(x\) 		&\(x^2\) 			& \(\so{\dots}\)			& \(x^{N-2}\) 			& \(x^{N-1}\) 			&\(\phantom{n}x^N\phantom{n}\)\\
	\hline
	\(y\)	\phantom{\Big|}		&  \(xy\) 	      	& \(x^2y\)			& \(\so{\dots}\) 			&\(x^{N-2}y\) 			& \(x^{N-1}y\) 			& 
	\tikz \node(char)[draw,fill=white,
  shape=rounded rectangle,
  minimum width=1cm,
  inner sep=2pt,
  transform canvas={yshift=.5mm}]
  {$x^Ny$};
    \\
    \noalign{\hrule height 1pt}
	\(y^2\) \phantom{\Big|}		& \(xy^2\) 	      	& \(x^2y^2\) 		& \(\so{\dots}\) 			& \(x^{N-2}y^2\) 		& \(-\)				& \(-\)\\
	\hline
	\(\so{\vdots}\) \phantom{\Big|}	& \(\so{\vdots}\) & \(\so{\vdots}\)		& \(\so{\rb{$\ddots$}}\)	& \(-\)	    			&\(\so{\rb{$\ddots$}}\)	& \(\so{\vdots}\)\\
	\hline
	\(y^{N-2}\) \phantom{\Big|}	& \(xy^{N-2}\) 	& \(x^2y^{N-2}\)		& \(-\)				& \(\so{\rb{$\ddots$}}\)	& \(-\) 				& \(-\)\\
	\hline
	\(y^{N-1}\) \phantom{\Big|}	& \({xy^{N-1}}\)  & \(-\) 			& \(\so{\rb{$\ddots$}}\)	& \(-\) 				& \(-\) 				& \(-\)\\
	\hline
	\(y^N\)	\phantom{$\Big|$}	& 
	\tikz \node(char)[draw,fill=white,
  shape=rounded rectangle,
  minimum width=1cm,
  inner sep=2pt,
  transform canvas={yshift=.5mm}]
   {$xy^N$};
	& \(-\) 			& \(\so{\dots}\)			& \(-\) 				&\(-\) 				& \(-\)\\
	\hline
	\end{tabular}
	\caption{Monomial basis for serendipity space \(S^N\). The circled elements extend $P^N$. This is used to define a minimal basis for $P^\mathrm{recon}$, see \eqref{Eq_PReconMin}. The bold lines separate the basis elements required for the edges from the triangle-shaped remainder that ensures unisolvence of the moments.}
	\label{Table_SerendipitySpaceBasis}
\end{table}

First, let us consider the case \(N \geq 4\) for the reconstruction space in \eqref{Eq_PReconMin} with $|S^N| = \tfrac12(N+1)(N+2)+2 = \{8, 12, 17, 23, \dots \}_{N \geq 2}$ degrees of freedom. There must be $4N$ point values (see Section \ref{Chap_DOFs}). This leaves us with \(\tfrac{1}{2}(N-3)(N-2) = |P^{N-4}|\) moments, that we will thus take with respect to the basis of $P^{N-4}$. 
For \(N = 2, 3\) with \(|S^N| = 4N\), there is an exception. In order to include the cell average, $S^N$ needs to be extended. In both cases \(\{x^2y^2\}\) is a suitable choice for the additional basis element. 
This choice of moments corresponds to the triangle-like choice \eqref{eq:trianglelike} and coincides with the tensor-like method for \(N+1=3\). 
An overview of the number of degrees of freedom can be found in Table \ref{Table_GenAFMinDOFs}.\\
\begin{table}
\fontsize{10pt}{10pt}\selectfont
\centering
\begin{tabular}{c | c c c}
	\hline
	order \phantom{$\Big|$} & $|P^\mathrm{recon}_{\min}|$ & \#PVs & \#Moments \\
	\!\!\!$N+1$ & $\displaystyle \frac{(N+1)(N+2)}{2}+2$ &  $4N$ & \(\displaystyle\max\left\{1, \frac{1}{2}(N-3)(N-2)\right\}\) \\
	\hline
	\hline
	3 \phantom{$\Big|$}& 8 \,\,\,\,\, 
	\tikz \node(char)[draw,fill=white,
  shape=rectangle,
  minimum width=.5cm,
  inner sep=2pt,
  transform canvas={yshift=.6mm}]
   {\textcolor{red}{\raisebox{.6mm}{+1}}};
   & 8   &  
	\tikz \node(char)[draw,fill=white,
  shape=rectangle,
  minimum width=.5cm,
  inner sep=2pt,
  transform canvas={yshift=.6mm}]
   {\textcolor{red}{\raisebox{.6mm}{1}}};
   \\
	4 \phantom{$\Big|$}& 12 \,\,\,\,\, 
	\tikz \node(char)[draw,fill=white,
  shape=rectangle,
  minimum width=.5cm,
  inner sep=2pt,
  transform canvas={yshift=.6mm}]
   {\textcolor{red}{\raisebox{.6mm}{+1}}};
   & 12 &  
	\tikz \node(char)[draw,fill=white,
  shape=rectangle,
  minimum width=.5cm,
  inner sep=2pt,
  transform canvas={yshift=.6mm}]
   {\textcolor{red}{\raisebox{.6mm}{1}}};
   \\
	5 \phantom{$\Big|$}& 17 & 16 & 1\\
	6 \phantom{$\Big|$}& 23 & 20 & 3\\
	7 \phantom{$\Big|$}& 30 & 24 & 6\\
	\hline
\end{tabular}
\caption{Number of degrees of freedom for GenAF$(M_{\symtri})$.}
\label{Table_GenAFMinDOFs} 
\end{table}

As mentioned above, the choice \(P^\mathrm{recon}=P^N\) (for $N\geq 5$) turned out not to be sufficient for unisolvence of the method on Cartesian cells, which was verified using \textsc{mathematica} \cite{Mathematica} for \(N = 5,6\), where we chose \(I_{\mathcal M} = I_{\mathcal M,{\symtri}} \setminus \{(N-4,0), (0,N-4)\}\). 
This is because two basis elements of the highest power \(N\) in \(x\) and \(y\) each are needed to be able to uniquely define the polynomials up to degree \(N\) on the two horizontal and vertical edges, respectively.
For \(N = 2,3,4\) it was also necessary to extend \(P^N\) in a suitable manner due to our requirements specified at the beginning of Section \ref{Chap_GenAF2DCartesian}. For $N=2,3$, we chose \(S^N\oplus \{x^2y^2\}\) as defined in \eqref{Eq_PReconMin} which we deemed a natural and symmetric choice lying within \(P^{N,N}\). For \(N=4\), we tried two versions: First, \(P^N \oplus \{x^3y^2, x^2y^3\}\), which we found insufficient to achieve unisolvence and second, we tested \(S^N\) as in \eqref{Eq_PReconMin}. 
Now, consider the minimal choice of degrees of freedom with the basis elements spanning \(P^\mathrm{recon}_{\min}\). Then, one can think of associating the nodal degrees of freedom (nodes and edge points) on the horizontal edges with the basis elements \(1, x, \dots, x^N\) and \(y, xy, \dots, x^Ny\) on an edge each, and respectively, \(1, y, \dots, y^N\) and \(x, xy, \dots, xy^N\) with the points on the vertical edges. The moments can then be associated with the basis elements in \eqref{eq:trianglelike}, which form a triangle in Table \ref{Table_SerendipitySpaceBasis} (indicated by the bold lines). 
As shown with the help of \textsc{mathematica} \cite{Mathematica} for up to order \(7\), the basis given in \eqref{Eq_PReconMin} for \(P^\mathrm{recon}_{\min}\) is unisolvent.\\

Referring to the choice of moments the following notation is introduced:
\begin{notation}
\label{Not_GenAFMMin}
	The generalized Active Flux method on two-dimensional Cartesian grids with \emph{minimal number of degrees of freedom} using $P^\mathrm{recon}_{\min}$ given in \eqref{Eq_PReconMin} shall be denoted with \emph{GenAF$(M_{\symtri})$}.
\end{notation}

\subsubsection{A hybrid finite element}
Next, the hybrid finite element as defined in \cite[Definition 3.1]{AB2023FEFV} is extended to a Cartesian element in two spatial dimensions. 
We also include the interpolation operator and the reconstruction and follow \cite[Definitions 5.2, 5.5, 5.7, 5.11, Proposition 5.12]{ErnGuermond2021} and \cite{AB2023FEFV}:
\begin{definition}\label{Def_FE}
Let $d=2$ and $N \geq 2$. Given an index set \(I_\mathrm{dof} := I_{\mathcal P} \cup I_{\mathcal M}\) with a set of points \(\{\mathbf x_p\}_{p \in I_{\mathcal P}}\), 
here \(\{\mathbf x_p = (x_p,y_p) | x_p, y_p \in \{\pm \tfrac{1}{2}\} \text{ or } (x_p \in \{\xi_a\}_a , y_p \in \{\pm \tfrac{1}{2}\}) \text{ or } (x_p \in \{\pm \tfrac{1}{2}\} , y_p \in \{\xi_a\}_a) \text{ with } \xi_a \in (-\tfrac{1}{2}, \tfrac{1}{2}) \quad \forall a = 0, \dots, N-2 \}\), and a (e.g. monomial) basis \(\{b_{k,l}\}_{k,l \in I_{\mathcal M}}\) for either
\begin{itemize}
	\item $P^{N-2,N-2}$ ($I_{\mathcal M}$ as in \eqref{eq:momentstensorlike}), or
	\item the polynomial vector space $P^{M}$ with $M = \max\{0, N-4\}$ ($I_{\mathcal M}$ as in \eqref{eq:trianglelike}).
\end{itemize}
Then, the finite element \((K, P, \Sigma)\) is 
\begin{itemize}
	\item \(K = C^\mathrm{ref} = \left[-\frac12, \frac12\right]\times\left[-\frac12, \frac12\right]\),
	\item $P$ given either by
	\begin{itemize}
		\item $P^\mathrm{recon} = P^{N,N}$, or
		\item $P^\mathrm{recon}=S^N \oplus \mathrm{span}\{x^2y^2\}$ for $N \leq 3$ and \(P^\mathrm{recon}=S^N\) for $N > 3$
	\end{itemize}
	and
	\item{$\Sigma =\{\sigma_i\}_{i \in I_\mathrm{dof}}$ the set of degrees of freedom with \(\sigma_i: P\rightarrow \mathbb R\)
	\begin{itemize}
		\item[]{\(\sigma_p(v) := v(\mathbf x_p) \quad \forall  p \in I_{\mathcal P}\)}
		\item[]{\(\sigma_{k,l}(v) := A^\mathrm{ref}_{k,l}\int_{K} b_{k,l}(\mathbf x) v(\mathbf x) \mathrm d \mathbf x \quad \forall (k,l) \in I_{\mathcal M}\)}
	\end{itemize}
	for all $v \in P$ and with the normalization factors \(A_{k,l}^\mathrm{ref}= (k+1)2^k(l+1)2^l\) with respect to $K$.
}
\end{itemize}
\end{definition}

\begin{figure}
\centering
	\begin{subfigure}[b]{0.45\textwidth}
	\centering
		\includegraphics[scale=0.425]{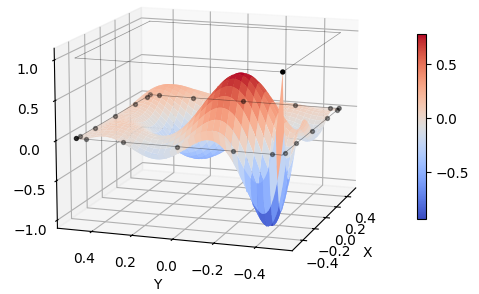}
		\caption{Node $(\mathcal N,0)$.}
		\label{Fig_GenAF7thShapeFuncNode} 
	\end{subfigure}
	\begin{subfigure}[b]{0.45\textwidth}
	\centering
		\includegraphics[scale=0.425]{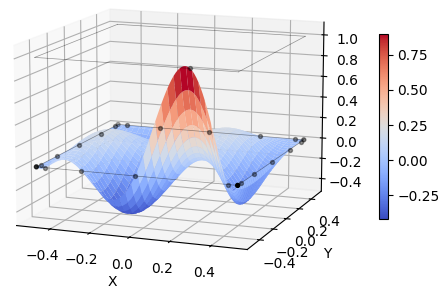}
		\caption{Edge point $(\mathcal E0,3)$.}
		\label{Fig_GenAF7thShapeFuncEdge} 
	\end{subfigure}
	\begin{subfigure}[b]{0.45\textwidth}
	\centering
		\includegraphics[scale=0.425]{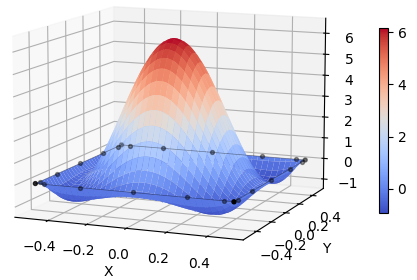}
		\caption{Cell average $(0,0)$.}
		\label{Fig_GenAF7thShapeFuncAv} 
	\end{subfigure}
	\begin{subfigure}[b]{0.45\textwidth}
	\centering
		\includegraphics[scale=0.425]{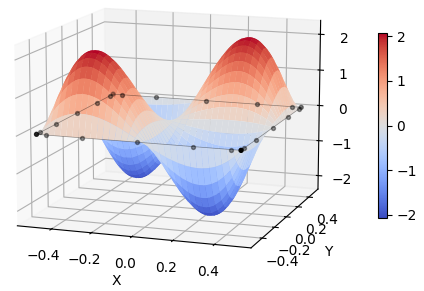}
		\caption{Moment $(1,1)$.}
		\label{Fig_GenAF7thShapeFuncMom} 
	\end{subfigure}
	\caption{Shape functions for GenAF$(M_{\symtri})$ of 7\textsuperscript{th} order.}
	\label{Fig_GenAF7thShapeFunctions} 
\end{figure}

We define the shape functions \(B_r: K\rightarrow \mathbb R\) such that
\begin{equation}
\sigma_r(B_s) = \delta_{rs} \quad \forall r,s \in I_\mathrm{dof}
\end{equation}
which form the canonical basis of $P$ with respect to the finite element. Figure \ref{Fig_GenAF7thShapeFunctions} shows examples for the shape functions of the GenAF$(M_{\symtri})$ method of order \(N+1 = 7\).

The interpolation operator corresponding to the finite element is defined to interpolate real-valued functions over $K$ in $P$. We consider the space $V:=L^1(K, \mathbb R)$.
With $P \subset V$, the elements in $\{\sigma_r\}_{r \in I_\mathrm{dof}}$ -- spanning the space of linear forms $\mathrm{hom}(P, \mathbb R)$ -- can naturally be extended to $\mathrm{hom}(V, \mathbb R)$. We denote them by $\{\tilde \sigma_r\}_{r \in I_\mathrm{dof}}$. 
Then, the interpolation operator $I_K: V \rightarrow P$ is defined as 
\begin{equation}\label{Eq_InterpolationOp}
	I_K(v)(\mathbf x) := \sum_{r \in I_\mathrm{dof}} \tilde \sigma_r(v) B_r(\mathbf x) \qquad \forall \mathbf x \in K.
\end{equation}
Finally, the reconstruction \(R: \mathbb R^{|I_\mathrm{dof}| \cdot s}\rightarrow P^s\),
\begin{equation}
	R\left((q_r)_{r \in I_\mathrm{dof}}\right)(\mathbf x) := \sum_{r \in I_\mathrm{dof}} q_{r} B_r(\mathbf x)\qquad \forall \mathbf x \in K
\end{equation}
maps the degrees of freedom $q_{r} \in I_\mathrm{dof}$ to the polynomial space $P^s$. 
This gives the reconstruction \eqref{Eq_DefRecon1}: 
\begin{align}
	q_\mathrm{recon}\vert_{C_{ij}} \circ \Phi^{-1}_{ij}&=R\left(\left(q_{ij}^{(r)}\right)_{r \in I_\mathrm{dof}}\right) \quad \in (P^\mathrm{recon})^s \\
	\Leftrightarrow \qquad  q_\mathrm{recon}|_{C_{ij}} (\mathbf x) &= \sum_{r \in I_\mathrm{dof}} q_{ij}^{(r)}B_r(\Phi_{ij}(\mathbf x)) .\label{Eq_DefRecon2}
\end{align}

\subsection{Update of point values}
\label{Chap_GenAF2DarbPV}
\begin{figure}
\centering
	\begin{subfigure}[b]{0.45\textwidth}
	\centering
		\includegraphics[scale=1.6]{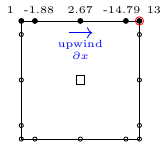}
		\caption{\((D_x^+)_{ij}^{(\mathcal N,2)}\).}
	\end{subfigure}
	\begin{subfigure}[b]{0.45\textwidth}
	\centering
		\includegraphics[scale=1.6]{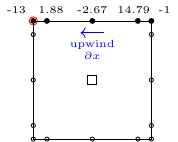}
		\caption{\((D_x^-)_{ij}^{(\mathcal N,2)}\).}
	\end{subfigure}

	\begin{subfigure}[b]{0.45\textwidth}
	\centering
		\includegraphics[scale=1.6]{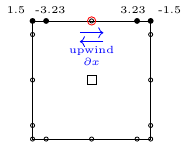}
		\caption{\((D_x)_{ij}^{(\mathcal E2,1)}\).}
	\end{subfigure}
	\begin{subfigure}[b]{0.45\textwidth}
	\centering
		\includegraphics[scale=1.6]{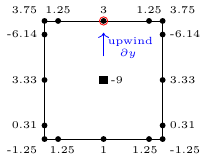}
		\caption{\((D_y^+)_{ij}^{(\mathcal E2,1)}\).}
	\end{subfigure}
	\begin{subfigure}[b]{0.45\textwidth}
	\centering
		\includegraphics[scale=1.6]{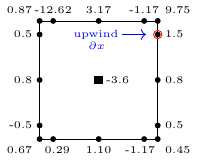}
		\caption{\((D_x^+)_{ij}^{(\mathcal E1,2)}\).}
	\end{subfigure}
	\begin{subfigure}[b]{0.45\textwidth}
	\centering
		\includegraphics[scale=1.6]{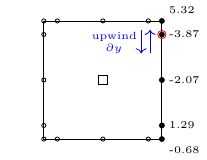}
		\caption{\((D_y)_{ij}^{(\mathcal E1,2)}\).}
	\end{subfigure}
	\caption{Example of finite difference formulas for GenAF$(M_{\symtri})$ of fifth order at locations of different degrees of freedom, marked by the additional circle (see \eqref{Eq_GenAFPVFDx}, \eqref{Eq_GenAFPVFDy}; coefficients rounded to two decimal points). Upwind direction: \enquote{\textcolor{blue}{$\rightarrow$}}.}
	\label{Fig_GenAF5PVFD}
\end{figure}

For the update of the point values it is possible to consider a non-conservative formulation of \eqref{Eq_SysCL}, which also allows to include stabilizing upwinding. In spirit we follow \cite{AB2023ExtensionAF, AB2023FEFV, ABK2025}. 
The quasi-linear form for the system of conservation laws \eqref{Eq_SysCL} is
\begin{equation}\label{Eq_QuasiLinCL}
	\partial_t q(t,\mathbf x)  +  Df^x(q(t,\mathbf x)) \partial_x q(t,\mathbf x) +  Df^y(q(t,\mathbf x) ) \partial_y q(t,\mathbf x)  = 0,
\end{equation}
which is considered at $\mathbf x_{ij}^{(p)}$ to derive the point value updates. 
Here, the Jacobians of the fluxes $f^x, f^y$ are given, since they can be calculated directly, whereas the derivatives of $q$ with respect to $x$ and $y$ have to be approximated. 
We intend to use different derivative approximations $(D_x)_{ij}^{(p)} \simeq \partial_x, (D_y)_{ij}^{(p)} \simeq \partial_y$ for each upwind direction. To this end, the diagonalized Jacobians 
\begin{align}
D f^x(q_{ij}^{(p)})&=T^x\mathrm{diag}(\lambda_{x,1},\dots, \lambda_{x,s})(T^x)^{-1},\\
D f^y(q_{ij}^{(p)})&=T^y\mathrm{diag}(\lambda_{y,1},\dots, \lambda_{y,s})(T^y)^{-1}
\end{align}
are split according to the positive and negative wave speeds to incorporate upwinding:
\begin{align*}
&D f^x(q_{ij}^{(p)})^{+} := T^x\mathrm{diag}(\lambda_{x,1}^{+},\dots , \lambda_{x,s}^{+})(T^x)^{-1},\\
&D f^x(q_{ij}^{(p)})^{-} := T^x\mathrm{diag}(\lambda_{x,1}^{-},\dots , \lambda_{x,s}^{-})(T^x)^{-1},\\
&D f^y(q_{ij}^{(p)})^{+} := T^y\mathrm{diag}(\lambda_{y,1}^{+},\dots , \lambda_{y,s}^{+})(T^y)^{-1},\\
&D f^y(q_{ij}^{(p)})^{-} := T^y\mathrm{diag}(\lambda_{y,1}^{-},\dots , \lambda_{y,s}^{-})(T^y)^{-1},
\end{align*}
where
\begin{equation*}
\lambda_{*,i}^{+} := \max(0, \lambda_{*,i}),  \qquad \lambda_{*,i}^{-} := \min(0, \lambda_{*,i}),\quad * = x, y.
\end{equation*}
This yields 
\begin{equation}
\label{Eq_PVUpdateJS}
\begin{aligned}
\frac{\mathrm d}{\mathrm d t}q_{ij}^{(p)} (t) =  &- Df^x(q_{ij}^{(p)})^{+}(D^{+}_x)_{ij}^{(p)} q - Df^x(q_{ij}^{(p)})^{-}(D^{-}_x)_{ij}^{({p})} q \\
&  -  Df^y(q_{ij}^{(p)})^{+}(D^{+}_y)_{ij}^{(p)} q -Df^y(q_{ij}^{(p)})^{-}(D^{-}_y)_{ij}^{(p)} q.
\end{aligned}
\end{equation}
The upwinded finite difference formulas $(D^{\pm}_x)_{ij}^{(p)}$, $(D^{\pm}_y)_{ij}^{(p)}$ are derived with the help of the reconstructions \(q_\mathrm{recon}|_{C_{ij}}\), which are differentiable on each cell \(C_{ij}\) (see also \cite{ABK2025, BKKL2025} for third order):
\begin{align}
(D^{\pm}_x)_{ij}^{(p)} q := \partial_x q_\mathrm{recon}\big\vert_{C^\mathrm{upw}}(\mathbf x)\big\vert_{\mathbf x_{ij}^{(p)}},\label{Eq_GenAFPVFDx}\\
(D^{\pm}_y)_{ij}^{(p)} q := \partial_y q_\mathrm{recon}\big\vert_{C^\mathrm{upw}}(\mathbf x)\big\vert_{\mathbf x_{ij}^{(p)}}.\label{Eq_GenAFPVFDy}
\end{align}
The upwind cell \(C^\mathrm{upw}\) is the cell adjacent to \(\mathbf x_{ij}^{(p)}\) from the corresponding \enquote{$\pm$}-direction for \(x\) or \(y\), i.e. at the node \(\mathbf x_{ij}^{(p)} = (x_{i+\frac{1}{2}}, y_{j+\frac{1}{2}})\)
\begin{align*}
&(D^{+}_x)_{ij}^{(p)} q:= \partial_x q_\mathrm{recon}|_{C_{ij}}(x, y_{j+\frac{1}{2}})|_{x=x_{i+\frac{1}{2}}}, \\
&(D^{-}_x)_{ij}^{(p)} q:= \partial_x q_\mathrm{recon}|_{C_{i+1,j}}(x, y_{j+\frac{1}{2}})|_{x=x_{i-\frac{1}{2}}},\\
&(D^{+}_y)_{ij}^{(p)} q:= \partial_y q_\mathrm{recon}|_{C_{ij}}(x_{i+\frac{1}{2}}, y)|_{y=y_{j+\frac{1}{2}}}, \\
&(D^{-}_y)_{ij}^{(p)} q:= \partial_y q_\mathrm{recon}|_{C_{i,j+1}}(x_{i+\frac{1}{2}}, y)|_{y=y_{j-\frac{1}{2}}},
\end{align*}
and at the horizontal edge points \((x_{i+\xi_{a}}, y_{j+\frac{1}{2}})\)
\begin{align*}
&(D_x)_{ij}^{(p)} q := (D^{+}_x)_{ij}^{(p)} q = (D^{-}_x)_{ij}^{(p)} = \partial_x q_\mathrm{recon}|_{C_{ij}}(x, y_{j+\frac{1}{2}})|_{x=x_{i+\xi_{a}}},\\
&(D^{+}_y)_{ij}^{(p)} q:= \partial_y q_\mathrm{recon}|_{C_{ij}}(x, y)|_{(x,y)=(x_{i+\xi_{a}},y_{j+\frac{1}{2}})},\\
&(D^{-}_y)_{ij}^{(p)} q:= \partial_y q_\mathrm{recon}|_{C_{i,j+1}}(x, y)|_{(x,y)=(x_{i+\xi_{a}},y_{j+\frac{1}{2}})}.
\end{align*}
The derivatives for the vertical edge points \((x_{i+\frac{1}{2}}, y_{j+\xi_{a}})\) are defined analogously. 
Considering definition \eqref{Eq_DefRecon2} for the reconstruction \(q_\mathrm{recon}|_{C_{ij}}\) the partial derivatives in $x$- and $y$-direction are given as
\begin{align}
\partial_x q_\mathrm{recon}|_{C_{ij}}(x,y) &= \sum_{r\in I_\mathrm{dof}}q_{ij}^{(r)}\nabla_\Phi B_r(\Phi_{ij}(x,y)) \cdot \partial_x \Phi_{ij}(x,y),\\
\partial_y q_\mathrm{recon}|_{C_{ij}}(x,y) &= \sum_{r\in I_\mathrm{dof}}q_{ij}^{(r)}\nabla_\Phi B_r(\Phi_{ij}(x,y)) \cdot \partial_y \Phi_{ij}(x,y).
\end{align}
An example for the thus derived finite difference formulas for \eqref{Eq_GenAFPVFDx} and \eqref{Eq_GenAFPVFDy} can be seen in Figure \ref{Fig_GenAF5PVFD}.

\subsection{Update of moments}
\label{Chap_GenAF2DarbMoments}
In the following, an update procedure is derived for the moments, extending \cite{AB2023ExtensionAF, AB2023FEFV, ABK2025}. In particular, a conservative update of the cell averages is obtained, which is key to ensure convergence to the weak solution of the conservation law. 
To this end, the weak formulation of \eqref{Eq_SysCL} 
\begin{equation} \label{Eq_WeakCL}
\begin{aligned}
	\frac{\mathrm d}{\mathrm d t} \int_{C}  v(\mathbf x) q (t, \mathbf x) \mathrm{d}C &+ \int_{\partial C}v(\mathbf x) \mathbf f\left(q(t, \mathbf x)\right) \cdot \mathbf n \mathrm{d}S  \\
	&-\int_{C}  \nabla v(\mathbf x) \cdot\mathbf f\left(q(t, \mathbf x)\right) \mathrm{d} \mathbf x  = 0
\end{aligned}
\end{equation}
on a domain $C \subset \mathbb R^d$, is considered. 
To find the update for the cell average and higher moments on cell $C_{ij}$, equation \eqref{Eq_WeakCL} with $C:= C_{ij}$ is multiplied by the normalization factor \eqref{Eq_GenAFNormalizationFactorMoments} and the test function \(v\) is set to \eqref{Eq_GenAFTestFunctionMoments}.
This yields the semi-discrete update formula for the moments of order $m=k+l$
\begin{equation}\label{Eq_WeakCLCell}
\begin{aligned}
	\frac{\mathrm d}{\mathrm dt} q_{ij}^{(k,l)}(t) =- A_{k,l}
	\bigg(
	& \int_{y_{j-\frac{1}{2}}}^{y_{j+\frac{1}{2}}} b_{k,l}\big(\Phi_{ij}(x_{i+\frac{1}{2}}, y)\big) f^x\big(q(t, x_{i+\frac{1}{2}}, y)\big) \\
	& \hspace{3em}	- b_{k,l}\big(\Phi_{ij}(x_{i-\frac{1}{2}}, y)\big)f^x\big(q(t, x_{i-\frac{1}{2}}, y)\big) \mathrm d y \\
	-& \int_{y_{j-\frac{1}{2}}}^{y_{j+\frac{1}{2}}} \int_{x_{i-\frac{1}{2}}}^{x_{i+\frac{1}{2}}}  \partial_x b_{k,l}\left(\Phi_{ij}(x, y)\right) f^x\left(q(t, x,y)\right) \mathrm d x \mathrm d y \\
	+& \int_{x_{i-\frac{1}{2}}}^ {x_{i+\frac{1}{2}}} b_{k,l}\big(\Phi_{ij}(x, y_{j+\frac{1}{2}})\big) f^y\big(q(t, x, y_{j+\frac{1}{2}})\big) \\
	& \hspace{3em}	- b_{k,l}\big(\Phi_{ij}(x, y_{j-\frac{1}{2}})\big) f^y\big(q(t, x, y_{j-\frac{1}{2}})\big) \mathrm d x& \\
	-& \int_{y_{j-\frac{1}{2}}}^{y_{j+\frac{1}{2}}} \int_{x_{i-\frac{1}{2}}}^{x_{i+\frac{1}{2}}}  \partial_y b_{k,l}\left(\Phi_{ij}(x, y)\right) f^y\left(q(t, x, y)\right) \mathrm d x \mathrm d y
	 \bigg), 
\end{aligned}
\end{equation}
where the integrals are approximated with the help of a numerical quadrature formula. Here, we choose either a Gauss-Lobatto or a Gauss-Legendre quadrature for the integrals. The necessary values at the quadrature points (both on the cell boundaries and within the cells) are obtained by evaluating the reconstruction \eqref{Eq_DefRecon2}. For now, the quadrature formulas are chosen to be sufficient for polynomials of at least degree $N + \max(k,l)$.\\
\noindent\emph{Remark:} The nodal degrees of freedom at the cell interfaces (nodes and edge points) could be included to calculate the quadratures. Yet, along the edge there are typically more point values than actually necessary for the quadrature alone. Here, we try to avoid this overintegration and the potentially larger computational effort. Furthermore, for the bulk integral, the reconstruction needs to be evaluated at quadrature points anyway as the pointwise degrees of freedom are located on cell interfaces only.

\subsection{Method definition}
\label{Chap_MethodDef}

Concluding, a definition is given of the generalized Active Flux method of arbitrarily high order on Cartesian grids in 2-d.
\begin{definition}\label{Def_GenAF}
Given $N \geq 2$, the \emph{generalized Active Flux method of order $N+1$} on 2-d Cartesian grids is the semi-discretization of \eqref{Eq_SysCL}\\
\begin{align*}
\left\lbrace
\begin{aligned}
	\frac{\mathrm d}{\mathrm dt} q_{ij}^{(k,l)}(t) =&- A_{k,l} \int_{\partial C_{ij}} b_{k,l}(\Phi_{ij}(\mathbf x)) \mathbf f(q(t,\mathbf x)) \cdot \mathbf n_{ij} \mathrm{d} S\\
									  &+ A_{k,l} \int_{C_{ij}} \nabla_{\mathbf x} b_{k,l}(\Phi_{ij}(\mathbf x)) \cdot \mathbf f(q(t,\mathbf x)) \mathrm d \mathbf x \quad \forall (k,l) \in I_{\mathcal M}\\
	\frac{\mathrm d}{\mathrm dt} q_{ij}^{(p)}(t) =& -F_p\left(\left(q_{i_p, j_p}^{(r)}\right)_{(i_p, j_p) \in \sigma_{p}, r \in I_\mathrm{dof}}\right) \quad \forall p \in \widetilde{I_{\mathcal P}}
\end{aligned}
\right.
\end{align*}
with $I_\mathrm{dof}$, $I_{\mathcal M}$, $\{b_{k,l}\}$ and $A^\mathrm{ref}_{k,l}$ given as in Definition \ref{Def_FE}, \(A_{k,l}=\tfrac{A^\mathrm{ref}_{k,l}}{\Delta x \Delta y}\) and the outward pointing unit normal $\mathbf n_{ij}$ on $\partial C_{ij}$. 
The set \(\big\{\mathbf x_{ij}^{(p)}\big\}_{p \in \widetilde{I_{\mathcal P}}}\) contains point values belonging to $C_{ij}$ as it is given by \eqref{Eq_IndexPointsBelCell}. 
$F_p$ is a consistent approximation of $\partial_x f^x(q) + \partial_y f^y(q)$ at point $\mathbf x_{ij}^{(p)}$ where $\sigma_p := \big\{(i,j) \big| \mathbf x_{ij}^{(p)} \in \partial C_{ij}\big\}$ defines the index set of the cell neighbours of $\mathbf x_{ij}^{(p)}$. 
\end{definition}

Here, we choose $F_p$ as given in \eqref{Eq_PVUpdateJS} with finite difference formulas \eqref{Eq_GenAFPVFDx}, \eqref{Eq_GenAFPVFDy} derived in Section \ref{Chap_GenAF2DarbPV}.

\section{Stability analysis}
\label{Chap_StabilityAnalysis}
Stability results for the generalized Active Flux method introduced in the previous Section are presented for linear advection in two spatial dimensions.
First, the semi-discrete generalized Active Flux method is studied with the help of an eigenvalue spectrum analysis. We observe that the stability of the method depends on the location of the edge points, and this analysis is used to determine a stable distribution. Second, a generalized Active Flux method using a third-order Runge-Kutta method for the discretization in time is considered to determine a stability bound for the fully discretized system. 
The methods used for the stability analysis are for example described in \cite{Lomax2001}.

\subsection{Semi-discrete Active Flux method}
\label{Chap_StabAnaSemiDiscreteGenAF}
The generalized Active Flux method of order $N+1$ uses moments and point values at the cell interfaces as degrees of freedom. In the one-dimensional case, the choice of the points is straightforward at the cell interfaces. The two-dimensional case allows more flexibility: out of $N+1$ points \(\mathbf x_{ij}^{(p)}\) along each cell interface, two are located at the nodes, (see Section \ref{Chap_DOFs}), but the question how the $N-1$ edge points are distributed remains.

Some first numerical tests with a uniform distribution of the edge points implied that this choice does not lead to a stable method. In \cite{AB2023ExtensionAF} a similar finding was observed for the arbitrarily high order one-dimensional Active Flux method with additional point values, reviewed in Section \ref{Chap_AFaddPV1d}. There, the additional points set within the cells had to be moved closer to the cell interfaces to gain stability. 
Inspired by this result we would like to find a similar way to distribute the edge points for the generalized Active Flux method in 2-d. 
Here, the focus is on symmetric distributions of the points with respect to the edge midpoint. For a fourth and fifth order method this leaves one free parameter for the edge point distribution and two for a sixth and seventh order method: Considering the edge \([-\frac12, \frac12]\) in $x$- or $y$-direction the edge points are 
\begin{itemize}
\item \(\{-\xi,\xi\}\) for \(N+1=4\) or \(\{-\xi, 0,\xi\}\) for \(N+1=5\), with \(\xi \in (0,\frac12)\) and 
\item \(\{-\xi_1, -\xi_0, \xi_0, \xi_1\}\) for \(N+1=6\) or \(\{-\xi_1, -\xi_0, 0, \xi_0, \xi_1\}\) for \(N+1=7\) with \(\xi_0, \xi_1 \in (0,\frac12)\), \(\xi_0 <\xi_1\)
\end{itemize}
While it might still be viable to analyze many different setups for a fourth or fifth order method, for higher orders \(N+1\) with \(\lfloor \tfrac{N-1}{2}\rfloor\) free parameters for the edge point distribution this is not feasible. Hence, the following stability analysis is confined to a few ways to distribute the edge points that we consider to be a sensible choice. 
In particular we will try Gauss and Gauss-Lobatto points. Even with Gauss points, we still include point values at the nodes in the set of our degrees of freedom. Although this choice is inspired by quadrature points, this choice is generally independent of the actual quadrature, see also Section \ref{Chap_GenAF2DarbMoments}. Also the number of point values is fixed with respect to the order of the method (see Section \ref{Chap_DOFs}).\\

Next, the eigenvalue spectrum of the semi-discrete generalized Active Flux method is analyzed for scalar linear advection in two spatial dimensions
\begin{equation}\label{Eq_linearAdvection}
\partial_t q(t, \mathbf x) + \nabla \cdot ( \mathbf a q(t, \mathbf x) )= 0
\end{equation}
with a constant advection speed \(\mathbf a = (a_x, a_y):= (\cos \theta, \sin \theta)^T\), \(\theta = [0, 2 \pi)\) and periodic boundary conditions. 
Thereto, the vector \(\mathbf q = (q_i)_{i \in \mathbf I_\mathrm{dofs}}\) with \(\mathbf I_\mathrm{dofs} = \bigcup_{i,j} (\widetilde {I_\mathrm{dof}})_{ij}\) assembling all degrees of freedom of the discretized domain \(\Omega_h\) is defined and the semi-discrete generalized Active Flux method for \eqref{Eq_linearAdvection} is re-written as the linear ODE system
\begin{equation}\label{Eq_ODEGenAF}
	\frac{\mathrm d}{\mathrm d t} \mathbf q(t) = \mathrm A \mathbf q(t),
\end{equation}
where $\mathrm A$ assembles the right hand side updates of the method.
A necessary condition for stability of the linear ODE system \eqref{Eq_ODEGenAF} is 
\begin{equation}
	\mathrm{Re}(\lambda_i) \leq 0 \qquad \forall \lambda_i \in \sigma(\mathrm A). 
\end{equation}
This condition becomes sufficient, if the algebraic multiplicity of all $\lambda_i$ with \(\mathrm{Re}(\lambda_i)=0\) equals their geometric multiplicity.\\

In the following, the GenAF$(M_{\symtri})$ method of orders \(N+1 = 3, \dots, 7\) is studied with Gauss, Gauss-Lobatto and uniformly distributed points on the edges. To cover a range of different advection directions, \(\theta \in [0,\tfrac{\pi}{2}]\) is chosen in increments of $\tfrac{1}{32}$-th, i.e. \(\theta = 0, \tfrac{\pi}{64} , \dots, \tfrac{\pi}{2}\), for the Gauss points and $\tfrac{1}{4}$-th for Gauss-Lobatto and uniformly distributed points. The restriction of \(\theta\) to \([0,\tfrac{\pi}{2}]\) is possible due to symmetry properties of the method. The setup for the grid sizes is $N_x=N_y \in \{ 3,5,10\}$.
The results use a Gauss integration sufficient for the calculation of the moments for linear advection, i.e. for moments of order $m$ it is chosen to be sufficient for polynomials of at least degree $N + \max(k,l)$.\\
\noindent\emph{Remark:} Some further tests with a higher accuracy or different type of quadrature yielded comparable results.\\

\begin{figure}
\centering
	\begin{subfigure}[b]{0.32\textwidth}
	\centering
		\includegraphics[scale=0.35]{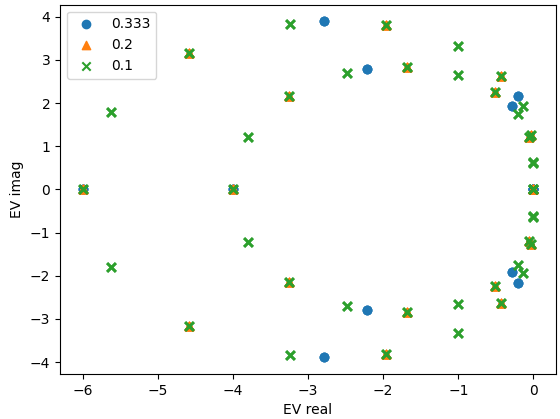}
		\caption{\(N+1 = 3\).}
	\end{subfigure}
	\begin{subfigure}[b]{0.32\textwidth}
	\centering
		\includegraphics[scale=0.35]{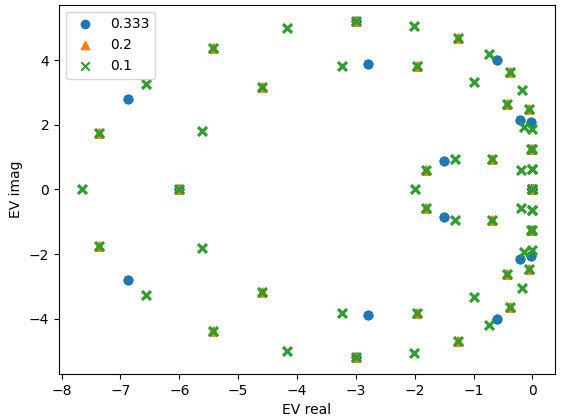}
		\caption{\(N+1 = 4\).}
	\end{subfigure}
	\begin{subfigure}[b]{0.32\textwidth}
	\centering
		\includegraphics[scale=0.35]{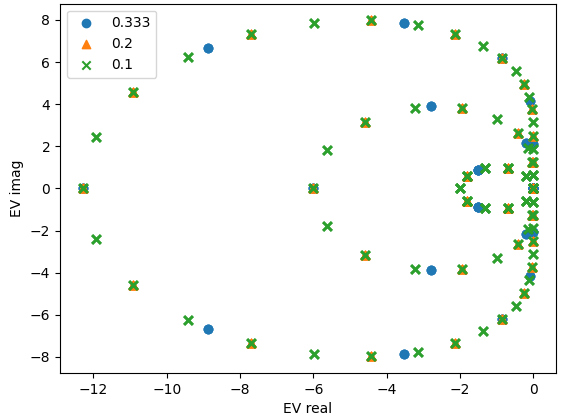}
		\caption{\(N+1 = 5\).}
	\end{subfigure}
	\begin{subfigure}[b]{0.32\textwidth}
	\centering
		\includegraphics[scale=0.35]{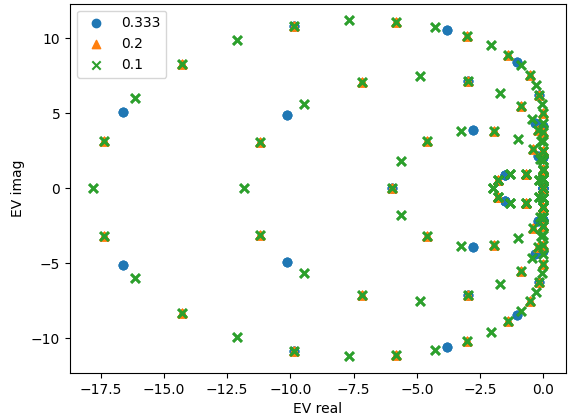}
		\caption{\(N+1 = 6\).}
	\end{subfigure}
	\begin{subfigure}[b]{0.32\textwidth}
	\centering
		\includegraphics[scale=0.35]{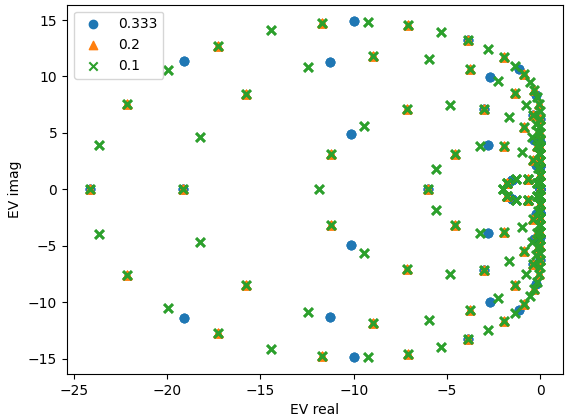}
		\caption{\(N+1 = 7\).}
	\end{subfigure}
	\caption{Scaled \(\lambda_i \in \sigma(\mathrm A)\), i.e. \(h\lambda_i\), for GenAF$(M_{\symtri})$ of order $N+1$ with Gauss edge point distribution. \(\theta = 0\) and \(h = \Delta x = \Delta y = 0.1, 0.2, \tfrac{1}{3}\). }
	\label{Fig_GenAFSemidiscreteStability}
\end{figure}

For the setups described above, the stability of \eqref{Eq_ODEGenAF} is analyzed numerically due to the size of the update matrices $\mathrm A$ (e.g. for seventh order, $10\times10$ grid: $1700\times1700$).\footnote{The matrix is assembled using an implementation of the method in C++ including the package \textsc{Eigen} \cite{cppeigenweb}. The actual analysis is done in Python using the \textsc{Numpy} package \cite{NumpyHarris2020} (v2.3) with its linear algebra library (\textsc{linalg}). For plotting we rely on the \textsc{Matplotlib} \cite{MatplotlibHunter2007} (v3.10).}
We obtain the following results: only the method with Gauss edge point distribution is stable for all tested orders \(N+1 = 3, \dots, 7\) and corresponding setups. From the computer-aided analysis one finds that \(\mathrm{Re}(\lambda_i) < \epsilon\), i.e. the eigenvalues have a non-postive real part up to a tolerance $\epsilon$, and all \(\lambda_i\) with \(|\mathrm{Re}(\lambda_i)| < \epsilon\) are at least semisimple for all $i$ with tolerance \(\epsilon = 5 \cdot 10^{-13}\) for $N+1 = 3,4,5$, \(\epsilon = 1 \cdot 10^{-12}\) for $N+1=6$ and \(\epsilon = 5 \cdot 10^{-12}\) for $N+1=7$ for the Gauss edge point distribution. In Figure \ref{Fig_GenAFSemidiscreteStability} an example is shown for \(N+1 = 3, \dots, 7\) with \(\theta = 0\) and the three grid sizes. It can also be seen that the eigenvalues scale with the grid size. 

The methods with uniform and with Gauss-Lobatto edge point distribution result in an unstable system \eqref{Eq_ODEGenAF} for orders \(N+1 = 4, \dots, 7\). For \(N+1 = 5, 6, 7\) the tests show that there exist eigenvalues with positive real part for all tested \(\theta\). For \(N+1 = 4\), eigenvalues with positive real part are observed for \(\theta = 0, \tfrac{\pi}{2}\). For \(N+1 = 3\) the edge point distribution coincides with the Gauss distribution.\\

In this way, we are led to conclude that Gauss points are a suitable choice for the distribution of the edge points that yield a stable method and we will use this distribution further on. 
Despite the limitations of the analysis, which has been conducted for orders \(N+1 = 3, \dots, 7\) and linear advection only, we are optimistic that these results will carry over to higher orders. Furthermore, we apply the Gauss edge point distribution to linear and non-linear systems of conservation laws in our numerical examples and have not observed instabilities.\\

\noindent\emph{Remark:} The semi-discrete GenAF$(M_{\square})$ method of orders \(N+1 = 4 \text{ and } 5\) was also studied for Gauss, Gauss-Lobatto and uniformly distributed edge points. The setup used for the advection directions is \(\theta \in [0, \tfrac{\pi}{2}], \theta = 0, \tfrac{\pi}{16}, \dots, \tfrac{\pi}{2}\) and the grid sizes are \(N_x=N_y=3,5,10\).
Again, only the Gauss edge point distribution yields a stable system \eqref{Eq_ODEGenAF} for all setups. For Gauss-Lobatto and uniformly distributed points, positive eigenvalues are found for \(\theta = 0, \tfrac{\pi}{2}\).

\subsection{Fully discrete Active Flux method}

The semi-discrete formulation of the method (see Definition \ref{Def_GenAF}) allows to choose a time discretization. Similarly as in \cite{AB2023FEFV, AB2023ExtensionAF} we focus on a strong stability preserving Runge-Kutta method of order three (SSP-RK3).\\

\begin{figure}
\centering
	\begin{subfigure}[b]{0.32\textwidth}
	\centering
		\includegraphics[scale=0.35]{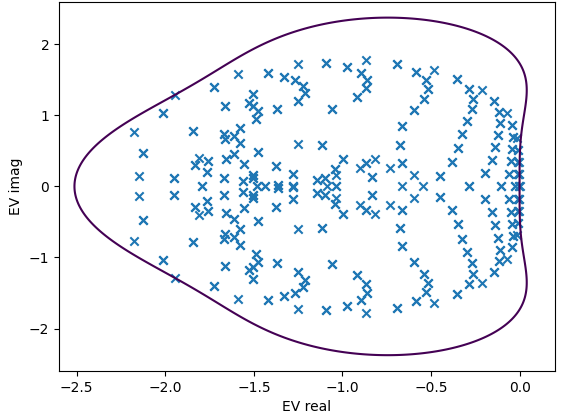}
		\caption{\(N+1 = 3\): \(\Delta t_\mathrm{max} \approx 0.038\).}
	\end{subfigure}
	\begin{subfigure}[b]{0.32\textwidth}
	\centering
		\includegraphics[scale=0.35]{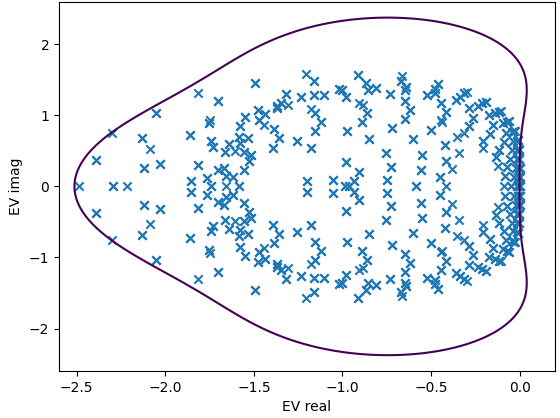}
		\caption{\(N+1 = 4\): \(\Delta t_\mathrm{max} \approx 0.029\).}
	\end{subfigure}
	\begin{subfigure}[b]{0.32\textwidth}
	\centering
		\includegraphics[scale=0.35]{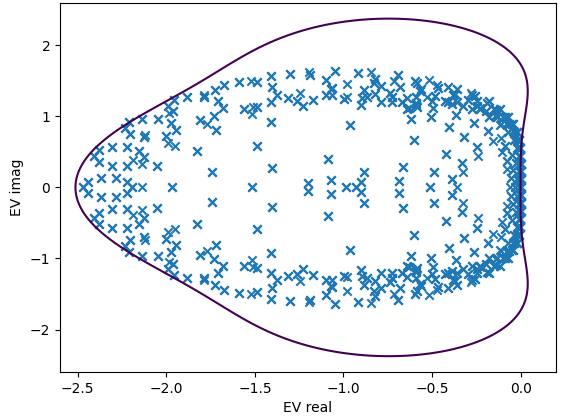}
		\caption{\(N+1 = 5\): \(\Delta t_\mathrm{max} \approx 0.025\).}
	\end{subfigure}
	\begin{subfigure}[b]{0.32\textwidth}
	\centering
		\includegraphics[scale=0.35]{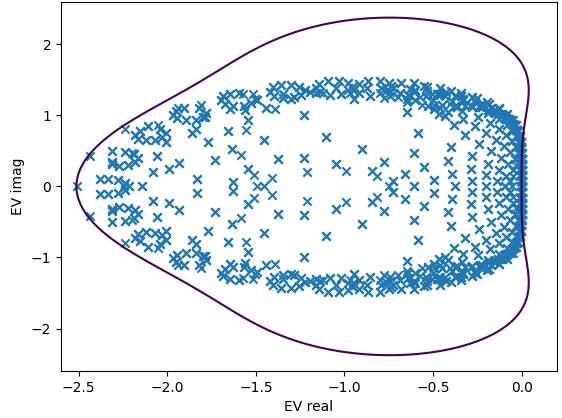}
		\caption{\(N+1 = 6\): \(\Delta t_\mathrm{max} \approx 0.018\).}
	\end{subfigure}
	\begin{subfigure}[b]{0.32\textwidth}
	\centering
		\includegraphics[scale=0.35]{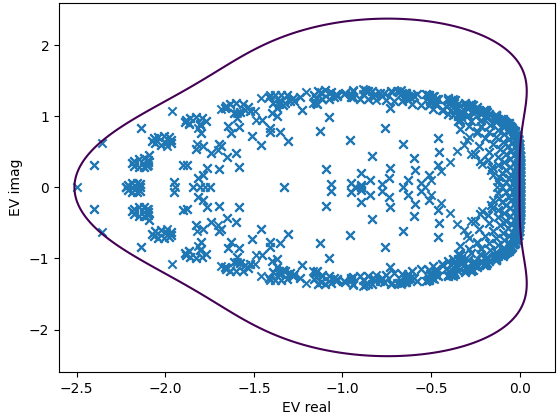}
		\caption{\(N+1 = 7\): \(\Delta t_\mathrm{max} \approx 0.013\).}
	\end{subfigure}
	\\
	\begin{subfigure}[b]{0.32\textwidth}
	\centering
		\includegraphics[scale=0.35]{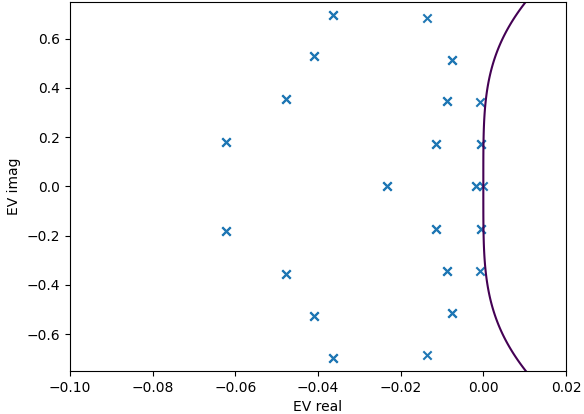}
		\caption{Close up for \(N+1 = 3\).}
	\end{subfigure}
	\begin{subfigure}[b]{0.32\textwidth}
	\centering
		\includegraphics[scale=0.35]{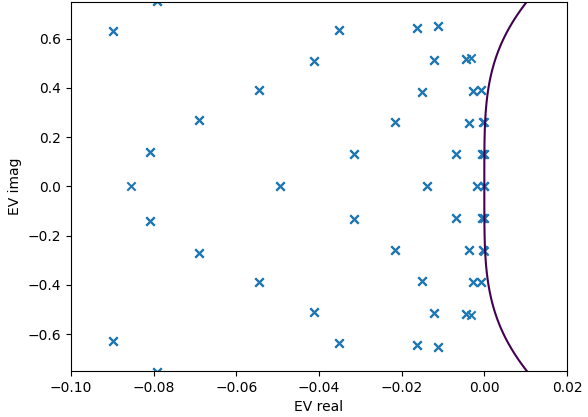}
		\caption{Close up for \(N+1 = 4\).}
	\end{subfigure}
	\begin{subfigure}[b]{0.32\textwidth}
	\centering
		\includegraphics[scale=0.35]{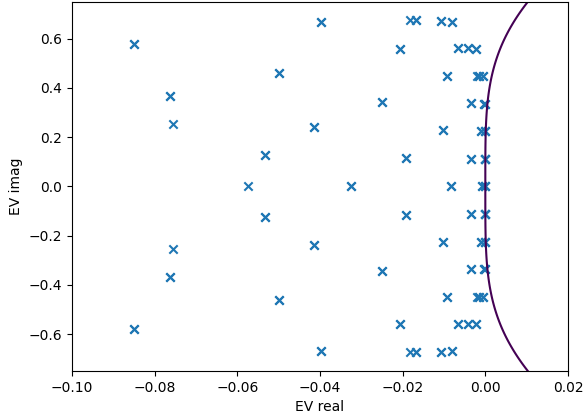}
		\caption{Close up for \(N+1 = 5\).}
	\end{subfigure}
	\begin{subfigure}[b]{0.32\textwidth}
	\centering
		\includegraphics[scale=0.35]{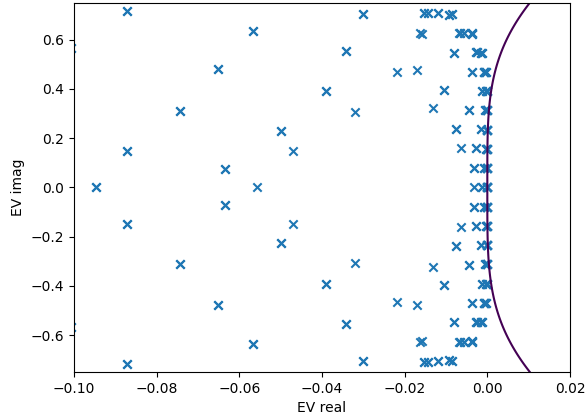}
		\caption{Close up for \(N+1 = 6\).}
	\end{subfigure}
	\begin{subfigure}[b]{0.32\textwidth}
	\centering
		\includegraphics[scale=0.35]{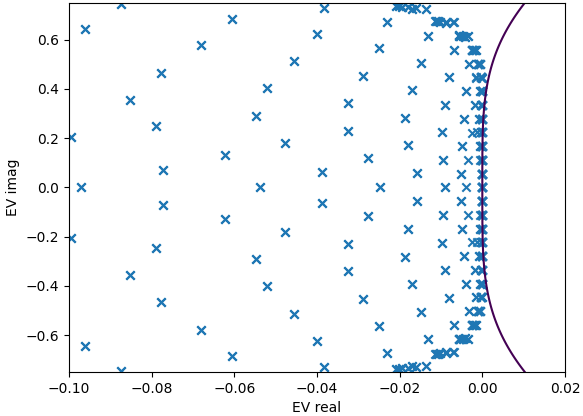}
		\caption{Close up for \(N+1 = 7\).}
	\end{subfigure}
	\caption{Stability domain $S_\text{RK3}$ (\textcolor{Purple}{\textbf{--}}) and scaled eigenvalues \(\lambda_i\Delta t_\text{max}\) (\textcolor{RoyalBlue}{$\times$}) for GenAF$(M_{\symtri})$ of order \(N+1\) in space for linear advection with \(\theta = \tfrac{\pi}{4}\), and \(\Delta x = \Delta y = 0.1\). Including close ups around origin.}
	\label{Fig_GenAFStabilityFullyDiscreteDeltatmax}
\end{figure}

Based on the approach in Section \ref{Chap_StabAnaSemiDiscreteGenAF} we analyze the stability of a fully discretized generalized Active Flux method with Gauss edge point distribution. This allows us to derive a CFL condition for linear advection \eqref{Eq_linearAdvection}. 
First, the ODE system \eqref{Eq_ODEGenAF} of the semi-discrete method for linear advection with periodic boundary conditions is diagonalized. This is possible, because the dimension of the sum of the eigenspaces for all eigenvalues of $\mathrm A$ is equal to $|\mathbf I_\mathrm{dofs}|$, which was computationally checked up to a tolerance $\epsilon = 5 \cdot 10^{-13}$ for $N+1 = 3,4,5$, \(\epsilon = 1 \cdot 10^{-12}\) for $N+1=6$ and \(\epsilon = 5 \cdot 10^{-12}\) for $N+1=7$. The diagonalized system for \(\hat{\mathbf q} = R^{-1} \mathbf q\) is given as
\begin{equation}\label{Eq_DiagODEGenAF}
\frac{\mathrm d}{\mathrm d t} \hat{\mathbf q}(t) = \Lambda \hat{\mathbf q}(t)
\end{equation}
with \(\mathrm A = R\Lambda R^{-1}\), \(\Lambda = \mathrm{diag}((\lambda_i)_{i \in {\mathbf I_\mathrm{dofs}}})\).
Applying the RK3 method with a time step \(\Delta t = t_{n+1} -t_n\) to \eqref{Eq_DiagODEGenAF} yields the fully discretized system
\begin{equation}\label{Eq_FullyDiscreted}
\mathbf{\hat q}_i^{n+1} = G(\lambda_i \Delta t) \mathbf{\hat q}_i^n \qquad \forall i \in {\mathbf I_\text{dofs}}
\end{equation}
with \(G(z):= 1+ z + \tfrac{1}{2} z^2 + \tfrac{1}{6} z^3\) and the stability domain \(S = \{z \in \mathbb C |  |G(z)| \leq 1\}\). 
From this it is possible to find a maximal time step \(\Delta t_\mathrm{max}(\theta)\), dependent on the advection direction, such that \(\lambda_i \Delta t_\mathrm{max}(\theta) \in S\) for all \(i \in {\mathbf I_\mathrm{dofs}}\). For our computer-based study, we increment the time step by \(10^{-4}\) for \(N+1=3,4,5\) and \(5\cdot 10^{-5}\) for \(N+1=6,7\) to approximate \(\Delta t_{\max}\). Figure \ref{Fig_GenAFStabilityFullyDiscreteDeltatmax} shows an example for the GenAF$(M_{\symtri})$ method of spatial order \(N+1 = 3, \dots, 7\) for \(N_x=N_y=10\). The scaled eigenvalues \(\lambda_i \Delta t_\mathrm{max} (\tfrac{\pi}{4})\) are plotted for the approximated maximal time step that still allows them to fit inside the stability domain for RK3 up to a tolerance $\epsilon$. 
\\

\begin{figure}
\centering
	\begin{subfigure}[b]{0.32\textwidth}
	\centering
		\includegraphics[scale=0.35]{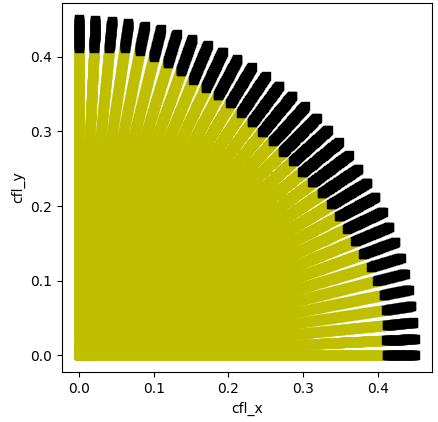}
		\caption{\(N+1 = 3\).}
	\end{subfigure}
	\begin{subfigure}[b]{0.32\textwidth}
	\centering
		\includegraphics[scale=0.35]{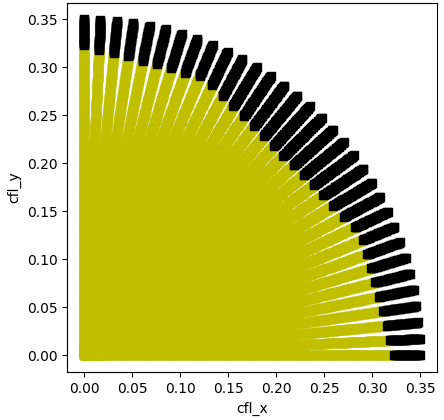}
		\caption{\(N+1 = 4\).}
	\end{subfigure}
	\begin{subfigure}[b]{0.32\textwidth}
	\centering
		\includegraphics[scale=0.35]{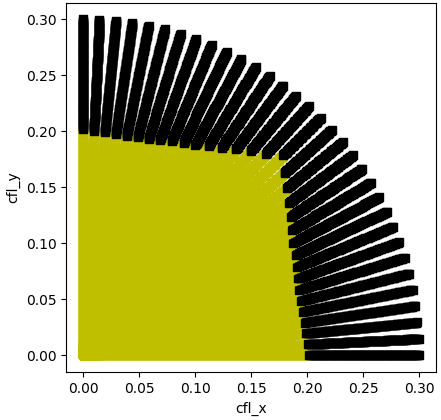}
		\caption{\(N+1 = 5\).}
	\end{subfigure}
	\begin{subfigure}[b]{0.32\textwidth}
	\centering
		\includegraphics[scale=0.35]{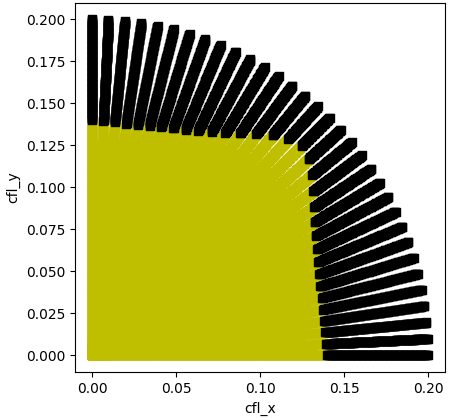}
		\caption{\(N+1 = 6\).}
	\end{subfigure}
	\begin{subfigure}[b]{0.32\textwidth}
	\centering
		\includegraphics[scale=0.35]{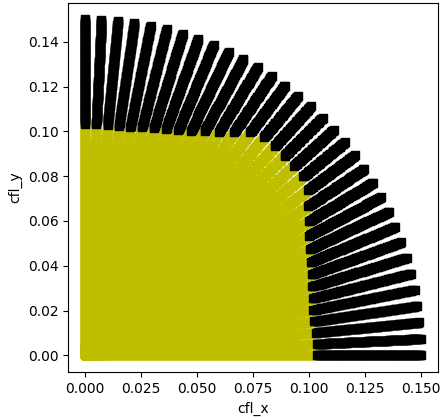}
	\caption{\(N+1 = 7\).}
	\end{subfigure}
	\caption{Stable regions (\textcolor{GreenYellow}{$\blacksquare$}) for GenAF$(M_{\symtri})$ of order \(N+1\) in space for linear advection. \(\Delta x = \Delta y = 0.1\) and RK3 time discretization. The yellow (\textcolor{GreenYellow}{$\blacksquare$}) and black ($\blacksquare$) domains depict the tested area.}
	\label{Fig_GenAFSSPRK3StableRegions}
\end{figure}

The stable regions for the fully discretized method of spatial order \(N+1 = 3, \dots,7\) are plotted as function of the CFLs \(\mathrm{cfl}_x = \frac{a_x \Delta t}{\Delta x}\) and \(\mathrm{cfl}_y = \frac{a_y \Delta t}{\Delta y}\) in $x$- and $y$-direction for \(h = \Delta x = \Delta y= 0.1\) in Figure \ref{Fig_GenAFSSPRK3StableRegions}. \\
One observes that the time step depends not only on the magnitude, but also the direction of the advection, and that the nature of this dependence changes as the order increases. This is studied below in more detail. In order to compare the CFL constraints for different orders of accuracy in a simple way, we consider the maximal CFL at \(a = a_x = a_y\)
\begin{equation}
C_\cfl:= \tfrac{|a|\Delta t_\text{max}}{h}
\end{equation}
given in Table \ref{Tab_CFLsGenAF2D} for the GenAF$(M_{\symtri})$ method discretized with SSP-RK3 in time, \(h = 0.1\).
\begin{table}
	\centering
	\begin{tabular}{c || c | c | c | c | c}
	order (\(N+1\)) & \(3\) & \(4\) & \(5\)& \(6\) & \(7\)\\
	\hline
	$C_\cfl$ (RK3) & \(0.27\) & \(0.20\) & \(0.17\)& \(0.12\) & \(0.088\)\\
	\end{tabular}
	\caption{\label{Tab_CFLsGenAF2D}\(C_\cfl\) for GenAF$(M_{\symtri})$ of order \(N+1\) in space for linear advection and RK3 discretization in time. (Two significant digits of $C_\cfl$ are shown for $h=0.1$.)}
\end{table}
For other grid sizes $h$ the results are comparable, because the eigenvalues $\lambda_i$ scale with $\tfrac{1}{h}$ (see Section \ref{Chap_StabAnaSemiDiscreteGenAF}).

The stability domain for the third-order method (with $M_{\symtri} = M_\square$) can approximately be described by
\begin{equation}\label{Eq_StableRegionCircle}
	S_\mathrm{circ} = \left\{(\mathrm{cfl}_x, \mathrm{cfl}_y) | \mathrm{cfl}_x^2 + \mathrm{cfl}_y^2 \leq r_\mathrm{cfl}^2 \right\}
\end{equation}
with \(r_\mathrm{cfl} := \sqrt{2} C_\cfl\) 
and gives the following bound for the time step size:
\begin{equation}\label{Eq_TimeStepBoundCircle}
	\Delta t \leq \frac{r_\cfl \min\{\Delta x, \Delta y\}}{\sqrt{a_x^2+a_y^2}}.
\end{equation}
For $N+1$=4 this is still a useful approximation while for higher spatial orders of the method the domain seems more restricted. Here, 
\begin{equation}\label{Eq_StableRegionSquare}
	S_{\mathrm{sq}} = \left\{(\mathrm{cfl}_x, \mathrm{cfl}_y) | \max\{\mathrm{cfl}_x,\mathrm{cfl}_y\} \leq C_\cfl \right\}
\end{equation}
could be considered which gives the bound
\begin{equation}
	\Delta t \leq \frac{C_\cfl \min\{\Delta x, \Delta y\}}{\max\{a_x, a_y\}}.
\end{equation}

\noindent\emph{Remark:} 
For $N+1=5$, GenAF$(M_{\symtri})$ almost yields a square. Looking more closely at the results in Figure \ref{Fig_GenAFSSPRK3StableRegions} this already shows for $N+1=4$, and for $N+1=6,7$ the corner starts to round off again. Comparing these results one can notice, that for GenAF$(M_{\symtri})$ with $N+1=4,5$ only point values at the cell interfaces are added as further degrees of freedom, seemingly favouring the maximal advection speed. Including higher moments for $N+1 \geq 6$ inside the cell decreases this effect.

\begin{figure}
	\begin{subfigure}[b]{0.45\textwidth}
	\centering
		\includegraphics[scale=0.35]{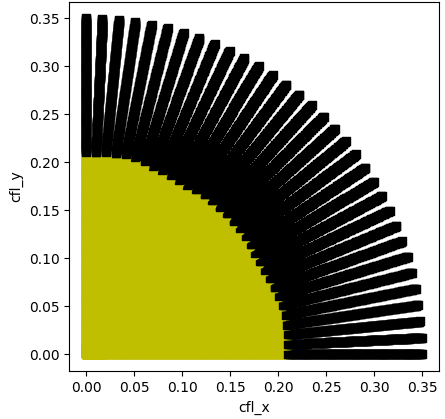}
		\caption{\(N+1 = 4\).}
	\end{subfigure}
	\begin{subfigure}[b]{0.45\textwidth}
	\centering
		\includegraphics[scale=0.35]{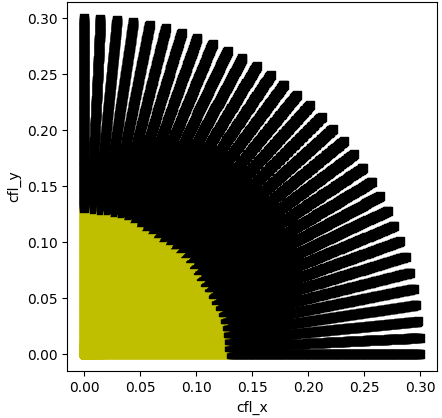}
		\caption{\(N+1 = 5\).}
	\end{subfigure}
	\caption{Stable regions (\textcolor{GreenYellow}{$\blacksquare$}) for GenAF$(M_{\square})$ of order \(N+1\) in space for linear advection. \(\Delta x = \Delta y = 0.1\) and RK3 time discretization. The yellow (\textcolor{GreenYellow}{$\blacksquare$}) and black ($\blacksquare$) domains depict the tested area.}
	\label{Fig_GenAFTensorSSPRK3StableRegions}
\end{figure}

\noindent\emph{Remark:} The GenAF$(M_{\square})$ method with Gauss edge point distribution and RK3 for the time discretization was analyzed for $N+1=4,5$. The stable regions \(S_{\text{GenAF}(M_\square)}\vert_{N+1=3,4,5}\) can be approximated by equation \eqref{Eq_StableRegionCircle}, see Figure \ref{Fig_GenAFTensorSSPRK3StableRegions}. The bounds \(r_\cfl\) are found to be $\sim 0.2$ ($C_\cfl \approx 0.14$) for $N+1=4$, and  $\sim 0.125$ ($C_\cfl \approx 0.089$) for $N+1=5$. Although \(C_\cfl\) is smaller than with $M_{\symtri}$, the stability domain with $M_\square$ seems to be less dependent on the advection direction. This fits with the observation above, considering the tensor-like added degrees of freedom.\\
Comparing the maximal CFL numbers to the one-dimensional case in \cite[Table 1 (Method A, RK3)]{AB2023FEFV} it can be observed that for the method in 2-d with $M_\square$ and RK3, the maximal CFL \(\approx r_\mathrm{cfl}\) for linear advection in $x$- or $y$-direction is approximately the same for the tested $N+1 = 3,4,5$.\\

In \cite{BKKL2025} a similar result for the generalized Active Flux method of third-order for linear acoustics in 2-d was indicated. A stability bound of $\mathrm{cfl} = \tfrac{\vert c \vert \Delta t}{\min \{\Delta x, \Delta y\}} < 0.28$ was derived using Fourier analysis.

\section{Numerical examples}
\label{Chap_NumEx}
Numerical results for the generalized Active Flux method GenAF$(M_{\symtri})$ for orders $3$ to $7$ on two-dimensional Cartesian grids with a SSP-RK3 time discretization are shown.
In the following we concentrate on smooth problems. More complex problems with e.g. discontinuities will need an appropriate limiting strategy, see e.g. \cite{DBK2025} for limiting of the third order Active Flux method. For higher orders this shall be part of future work.

\subsection{Linear advection}
\label{Chap_NumExLinAd}
\begin{figure}
	\centering
	\begin{subfigure}[b]{0.45\textwidth}
		\centering
		\includegraphics[scale=0.425]{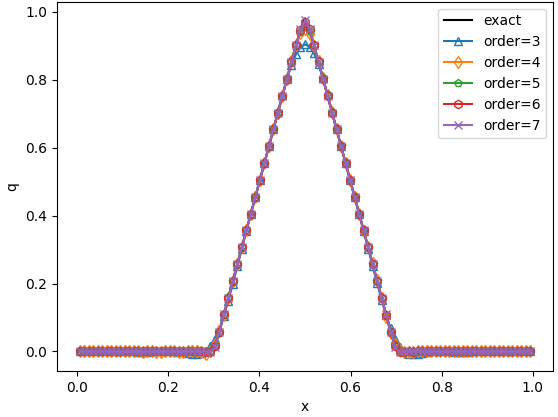}
	\end{subfigure}
	\begin{subfigure}[b]{0.45\textwidth}
		\centering
		\includegraphics[scale=0.425]{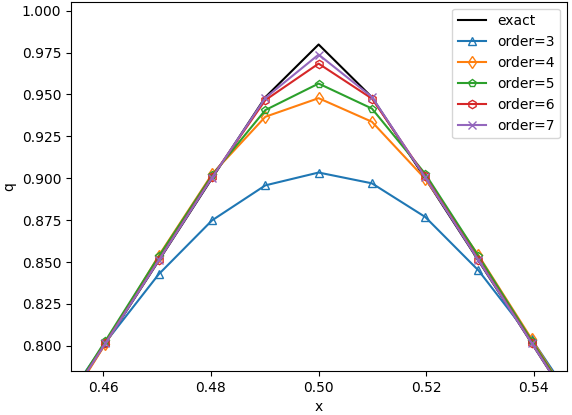}
	\end{subfigure}
	\caption{\label{Fig_NumExLinAdCone01} Example of diagonally advected cone. Cross section at $y=0.5$ of cell average on grid \(101\times101\) at \(t= 5\) for GenAF$(M_{\symtri})$ of order \(N+1=3, \dots, 7\) (left) and close up (right).}
\end{figure}

The first example considers linear advection \(\partial_t q(t, \mathbf x) + \partial_x q(t, \mathbf x) + \partial_y q(t, \mathbf x) = 0\) on \([0,1]^2\) with periodic boundary conditions where initial data in the shape of a cone of radius \(r_\mathrm{max}=0.2\)
\begin{align*}
	q(0, \mathbf x) =
	\left\lbrace
	\begin{aligned}
	&1- \frac{r}{r_\text{max}} \quad &&\text{for } r < r_\mathrm{max},\\
	&0 &&\text{otherwise}
	\end{aligned}
	\right.
\end{align*}
with \(r = \sqrt{(x-0.5)^2+(y-0.5)^2}\) is being advected diagonally through the domain. Figure \ref{Fig_NumExLinAdCone01} shows a cross section of the solution at \(t = 5\) on a grid with \(101 \times 101\) cells for spatial orders $N+1=3$ to $7$. A $\mathrm{cfl} = C_\cfl$ as shown in Table \ref{Tab_CFLsGenAF2D} (except for $N+1=7$: $\mathrm{cfl} = 0.085$) is used to calculate the time step. One observes that the apex of the cone is better approximated the higher the order of the method.
\\

\begin{figure}
	\centering
	\includegraphics[scale=0.6]{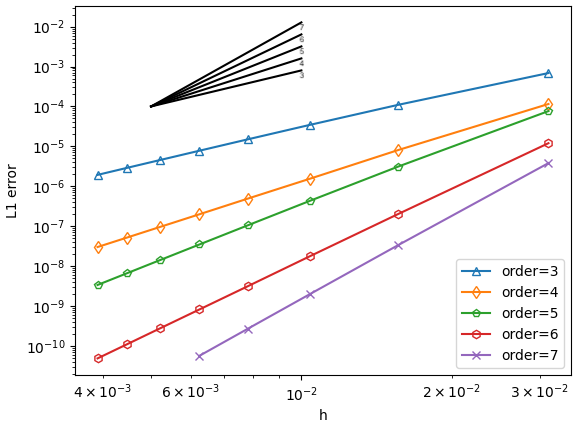}
	\caption{Convergence of the cell averages in $L^1$-error for linear advection for GenAF$(M_{\symtri})$ of order \(N+1=3, \dots, 7\) and RK3 time discretization.}
	\label{Fig_NumExLinAdConvergence} 
\end{figure}

The second example considers smooth initial data (in resemblance to \cite[Section 3.5.1]{AB2023ExtensionAF})
\begin{equation*}
q(0,\mathbf x) = 0.8+\exp{\left(-\left(\tfrac{x-0.5}{0.05}\right)^2 - \left(\tfrac{y-0.5}{0.05}\right)^2\right)}
\end{equation*}
that are used to show the convergence order of the method with \(N+1 = 3, \dots, 7\) at \(t=0.1\). To recover the spatial convergence order $N+1$ with a RK3 time discretization \(\Delta t^3 = \mathcal O(h^{N+1})\) is needed, \(h = \Delta x= \Delta y\). For the cell sizes \(\{h_i\}_i\) an adaptive \(\mathrm{cfl}(h_i) = \mathrm{cfl}(h_1) (\tfrac{h_i}{h_1})^{\frac{N-2}{3}}\) is used starting with $h_1 = \tfrac{1}{32}$ and the \(\mathrm{cfl}(h_1) = C_\cfl\) as above. In Figure \ref{Fig_NumExLinAdConvergence} and the corresponding Table \ref{Tab_ConvOrderGenAF2DMin} the convergence of the error of the cell averages in the $L^1$-norm is shown and a convergence order \(\mathcal O(h^{N+1})\) is observed.

\begin{table}
	\fontsize{11pt}{11pt}\selectfont
	\centering
	\begin{tabular}{l || c | c | c | c | c | c |}
	\(N+1\) \phantom{\Big|}& \multicolumn{2}{c|}{\(3\)} & \multicolumn{2}{c|}{\(4\)} & \multicolumn{2}{c|}{\(5\)}\\
	\hline
	
	$h$ \phantom{\big |}& $e_\mathrm{L1}$ & EOC & $e_\mathrm{L1}$ & EOC & $e_\mathrm{L1}$ & EOC\\
	\hline
	$0.03125$ 	& $6.87 \cdot 10^{-4}$ & - 	      & $1.15 \cdot 10^{-4}$ & - 	    & $7.65 \cdot 10^{-5}$ & - 	  	\\
	$0.015625$ 	& $1.10 \cdot 10^{-4}$ & $2.65$ & $8.06 \cdot 10^{-6}$ & $3.84$ & $3.10 \cdot 10^{-6}$ & $4.62$ \\
	$0.0104167$ 	& $3.46 \cdot 10^{-5}$ & $2.84$ & $1.55 \cdot 10^{-6}$ & $4.07$ & $4.33 \cdot 10^{-7}$ & $4.86$ \\
	$0.0078125$ 	& $1.50 \cdot 10^{-5}$ & $2.91$ & $4.89 \cdot 10^{-7}$ & $4.01$ & $1.05 \cdot 10^{-7}$ & $4.94$ \\
	$0.00625$ 	& $7.76 \cdot 10^{-6}$ & $2.95$ & $1.98 \cdot 10^{-7}$ & $4.05$ & $3.46 \cdot 10^{-8}$ & $4.96$ \\
	$0.00520833$ 	& $4.52 \cdot 10^{-6}$ & $2.96$ & $9.50 \cdot 10^{-8}$ & $4.03$ & $1.40 \cdot 10^{-8}$ & $4.97$ \\
	$0.00446429$ 	& $2.86 \cdot 10^{-6}$ & $2.97$ & $5.11 \cdot 10^{-8}$ & $4.03$ & $6.49 \cdot 10^{-9}$ & $4.98$ \\
	$0.00390625$ 	& $1.92 \cdot 10^{-6}$ & $2.98$ & $2.98 \cdot 10^{-8}$ & $4.02$ & $3.34 \cdot 10^{-9}$ & $4.98$ \\
	\hline
	\hline
	
	\(N+1\)\phantom{\Big|} & \multicolumn{2}{c|}{\(6\)} & \multicolumn{2}{c|}{\(7\)} & \multicolumn{2}{c}{}\\
	\cline{1-5}
	$h$ \phantom{\big |}& $e_\mathrm{L1}$ & EOC & $e_\mathrm{L1}$ & EOC & \multicolumn{2}{c}{}\\
	\cline{1-5}
	$0.03125$ 	&  $1.20 \cdot 10^{-5}$ & - 	     & $3.79 \cdot 10^{-6}$ & - & \multicolumn{2}{c}{}\\
	$0.015625$ 	&  $2.01 \cdot 10^{-7}$ & $5.90$ & $3.33 \cdot 10^{-8}$ & $6.83$ & \multicolumn{2}{c}{}\\
	$0.0104167$ 	&  $1.77 \cdot 10^{-8}$ & $5.99$ & $1.99 \cdot 10^{-9}$ & $6.95$ & \multicolumn{2}{c}{}\\
	$0.0078125$ 	&  $3.11 \cdot 10^{-9}$ & $6.05$ & $2.67 \cdot 10^{-10}$ & $6.98$ & \multicolumn{2}{c}{}\\
	$0.00625$ 	& $8.13 \cdot 10^{-10}$ & $6.01$ & $5.60 \cdot 10^{-11}$ & $7.00$ & \multicolumn{2}{c}{}\\
	$0.00520833$ 	& $2.72 \cdot 10^{-10}$ & $6.01$ & - 			& - & \multicolumn{2}{c}{}\\
	$0.00446429$ 	& $1.07 \cdot 10^{-10}$ & $6.03$ & - 			& - & \multicolumn{2}{c}{}\\
	$0.00390625$ 	& $4.81 \cdot 10^{-11}$ & $6.01$ & - 			& - & \multicolumn{2}{c}{}\\
	\end{tabular}
	\caption{\label{Tab_ConvOrderGenAF2DMin} $L^1$-error $e_\mathrm{L1}$ and corresponding experimental order of convergence (EOC) for GenAF$(M_{\symtri})$ of order \(N+1\) of the cell averages for linear advection and RK3 time discretization.}
\end{table}

\subsection{Acoustic equations}
\begin{figure}
	\centering
	\includegraphics[scale=0.4]{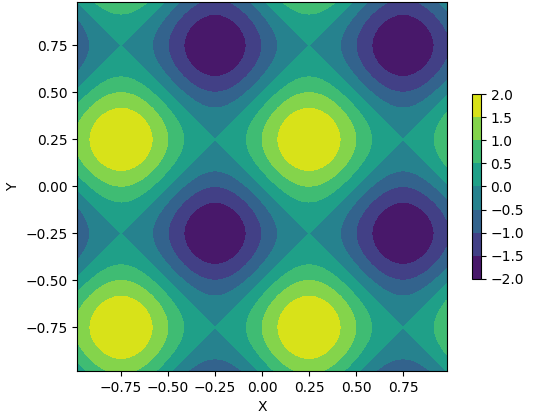}
	\caption{Example for acoustics with sine wave in pressure. Exact solution for the cell average $p^{(0,0)}$ of pressure on a \(60\times60\) grid at \(t= 5\). \emph{Remark:} The cell average values are approximated by numerical integration.}
	\label{Fig_NumExAcousticsPressureSinwaveExSoltend1}
\end{figure}

\begin{figure}
	\centering
	\begin{subfigure}[b]{0.32\textwidth}
		\centering
		\includegraphics[scale=0.33]{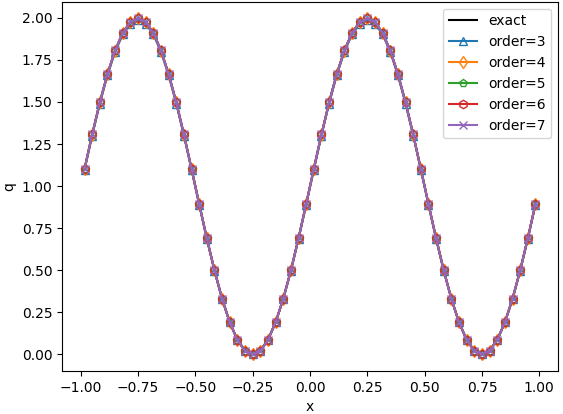}
	\end{subfigure}
	\begin{subfigure}[b]{0.32\textwidth}
		\centering
		\includegraphics[scale=0.33]{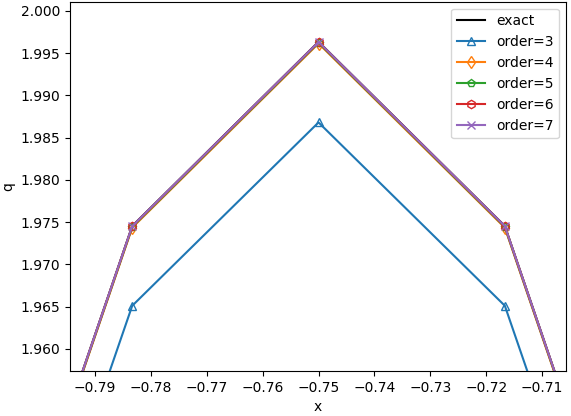}
	\end{subfigure}
	\begin{subfigure}[b]{0.32\textwidth}
		\centering
		\includegraphics[scale=0.33]{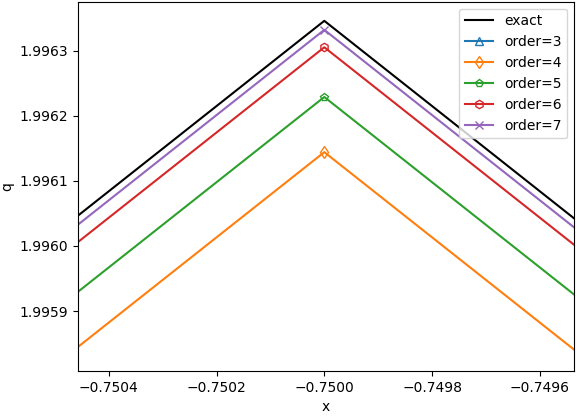}
	\end{subfigure}
	\caption{Example for acoustics with sine wave in pressure. Cross section at $y=0.25$ of $p^{(0,0)}$ on grid \(60\times60\) at \(t= 5\) for GenAF$(M_{\symtri})$ of order \(N+1=3, \dots, 7\) (left) and two close ups (middle, right) zooming in at the left maximum (observe the different scales of the $x$-axes).}
	\label{Fig_NumExAcousticsPressureSinwave} 
\end{figure}

The following example for the acoustic equations with \(c>0\)
\begin{align*}
	\partial_t p + c\nabla \cdot \mathbf v &= 0 \\
	\partial_t \mathbf v +  c\nabla p &= 0 
\end{align*}
on $[-1, 1]^2$ considers initial conditions with a sine wave in pressure and zero initial velocity as suggested in \cite{ER2013} from \cite{LMW2000} 
\begin{align*}
	p(0, \mathbf x) &= \frac{1}{c}(\sin(2\pi x)+\sin(2\pi y)),\\
	\mathbf v(0, \mathbf x) &= (0, 0)^{\mathrm T}
\end{align*}
with periodic boundary conditions, for which the exact solution is given as
\begin{align*}
	p(t, \mathbf x) &= \frac{1}{c}\cos(2\pi ct)(\sin(2\pi x)+\sin(2\pi y)),\\
	\mathbf v(t, \mathbf x) &= \frac{1}{c}\sin(2\pi ct)(\cos(2\pi x), \cos(2\pi y))^{\mathrm T}.
\end{align*}
We consider $c=1$. 
The solution at \(t= 5\) where the exact solution matches the initial condition is computed (see Figure \ref{Fig_NumExAcousticsPressureSinwaveExSoltend1}). The grid size is \(N_x = N_y = 60\) and the time step \(\Delta t = \mathrm{cfl} \tfrac{\min\{\Delta x, \Delta y\}}{c}\) is computed with $\mathrm{cfl} = C_\cfl$, a CFL bound obtained for linear advection, see Table \ref{Tab_CFLsGenAF2D} (except for $N+1=7$: $\mathrm{cfl} = 0.085$). As for linear advection, a better approximation of the exact solution is observed for higher orders of the generalized Active Flux method. In particular, a cross section through a maximum of the pressure sine wave is shown in Figure \ref{Fig_NumExAcousticsPressureSinwave}.

\subsection{Euler equations}

\begin{figure}
	\centering
	\includegraphics[scale=0.4]{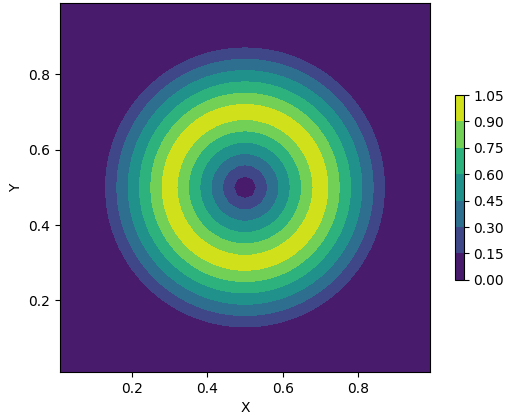}
	\caption{\label{Fig_NumExEulerGreshoCellAvMomExSoltend1} Gresho vortex for Euler equations. Exact solution for the norm of the cell average of the momentum $\Vert (\rho \mathbf v)^{(0,0)}\Vert_2$ on a \(51\times51\) grid at \(t= 1\). \emph{Remark:} The cell average values are approximated by numerical integration.}
\end{figure}

\begin{figure}
	\centering
	\begin{subfigure}[b]{0.45\textwidth}
		\centering
		\includegraphics[scale=0.425]{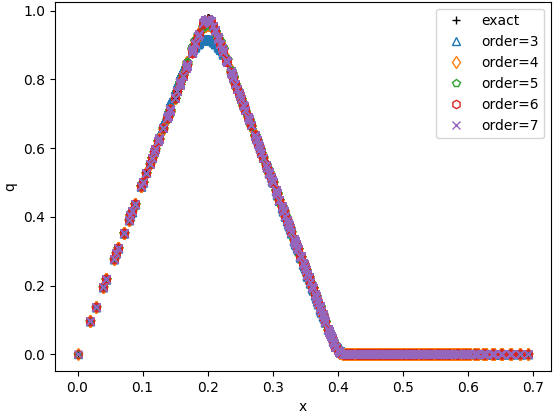}
	\end{subfigure}
	\begin{subfigure}[b]{0.45\textwidth}
		\centering
		\includegraphics[scale=0.425]{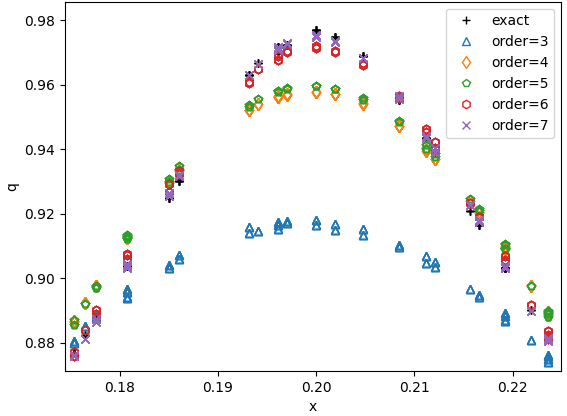}
	\end{subfigure}
	\caption{Radial plot of $\Vert (\rho \mathbf v)^{(0,0)}\Vert_2$ for Gresho vortex on \(51\times51\) grid at \(t= 1\) for GenAF$(M_{\symtri})$ of order \(N+1=3, \dots, 7\) (left) and close up (right).}
	\label{Fig_NumExEulerGreshoCellAvMomEx} 
\end{figure}

The generalized Active Flux method can also be applied to non-linear problems like the compressible Euler equations
\begin{align*}
	\partial_t \rho + \nabla \cdot (\rho \mathbf v) &= 0 \\
	\partial_t (\rho \mathbf v) +  \nabla \cdot (\rho \mathbf v \otimes \mathbf v + p \mathbb I) &= 0 \\
	\partial_t E + \nabla \cdot ((E+p)\mathbf v) &= 0
\end{align*}
In this paper an ideal polytropic gas is considered with \(E = \frac{p}{\gamma - 1} +\frac{1}{2}\rho \lVert \mathbf v \rVert^2_2 \) and \(\gamma =1.4\).\\
A Gresho vortex (see \cite{GC1990, BEKMR2017aa}) on $[0, 1]^2$
\begin{align*}
	\rho(0, \mathbf x) &= 1\\
	\mathbf v(0, \mathbf x) &= {\left\lbrace
	\begin{aligned}
	&\tfrac{5r}{r}(-y-0.5,x-0.5)^\mathrm{T} &\qquad  & r < 0.2\\
	&\tfrac{2-5r}{r}(-y-0.5,x-0.5)^\mathrm{T}& \qquad & 0.2 \leq r < 0.4\\
	&0& \qquad & \text{else}\\
	\end{aligned}
	\right.
	}\\ 
	p(0, \mathbf x) &=\left\lbrace
	\begin{aligned}
	 &\tfrac{1}{\gamma M^2}-\tfrac{1}{2}+ \tfrac{(5r)^2}{2} & \qquad &r < 0.2\\
	 &\tfrac{1}{\gamma M^2}-\tfrac{1}{2}+ 4\ln(5r)+4-20r+ \tfrac{(5r)^2}{2} & \qquad 0.2 &\leq r < 0.4\\
	 &\tfrac{1}{\gamma M^2}-\tfrac{1}{2}+ 4\ln2-2 & \qquad &\text{else}\\
	\end{aligned}
	\right.
\end{align*}
for $M=0.1$ with \(r = \sqrt{(x-0.5)^2+(y-0.5)^2}\) and periodic boundary conditions is computed at \(t = 1\) on a grid of \(51\times 51\) cells (see Figure \ref{Fig_NumExEulerGreshoCellAvMomExSoltend1}). The time step is calculated with the help of the CFL bound from Table \ref{Tab_CFLsGenAF2D} (except for $N+1=7$: $\mathrm{cfl} = 0.085$). 
Figure \ref{Fig_NumExEulerGreshoCellAvMomEx} shows the radial plot of the cell average for the norm of the momentum \(\Vert (\rho \mathbf v)^{(0,0)}\Vert_2\) for the generalized Active Flux method of order \(3, \dots, 7\). 

\section{Conclusion and outlook}
This paper presents a generalized Active Flux method on two-dimensional Cartesian grids of arbitrarily high order using higher moments in addition to the cell averages and point values at the cell interfaces. It is extending a semi-discrete hybrid finite element--finite volume version of Active Flux from 1-d. We focused on a method with a serendipity-like, hybrid finite element using a minimal number of degrees of freedom on the Cartesian cells to keep the computational cost minimal. An alternative is a tensor-like extension of the element from 1-d to 2-d. The point values at the cell interfaces include the nodes and edge points, such that a globally continuous reconstruction is obtained. Their update uses a Jacobian splitting of the non-conservative formulation of the system and finite difference formulas for the approximation of the derivatives of the solution. These are derived using the reconstruction of the solution on the cells. The moments are updated based on the weak formulation with the integrals approximated by quadrature formulas. 

The eigenvalue spectrum of the semi-discrete method was analyzed for linear advection up to order $7$. We find that the edge points at the cell interfaces cannot be distributed arbitrarily. We have identified, at least for orders up to 7, Gauss points as a suitable choice to achieve a stable setup of the semi-discrete method. For the time discretization we relied on a SSP-RK3 method and found CFL bounds of the method up to order $7$ for linear advection. Numerical examples confirm the theoretical convergence order and demonstrate how the high-order methods allow to better resolve features of the solution. 

Further studies to better understand the effect of suitable edge point distributions on the stability of the method will be interesting.
In future work, it will also be important to develop a suitable limiting for the method. To this end, a flux vector splitting for the point value update, which was recently presented in \cite{DBK2025} for third order, could be considered.
Further work shall also consider an extension to 3-d as well as applications to simulations of turbulent convection and wave propagation, which require very high orders of accuracy.

\section*{Acknowledgements}
\noindent WB, CK and LL acknowledge funding by the Deutsche Forschungsgemeinschaft (DFG, German Research Foundation) within \emph{SPP~2410 Hyperbolic Balance Laws in Fluid Mechanics: Complexity, Scales, Randomness (CoScaRa)}, project number 525941602.\\
PC acknowledges support from Department of Atomic Energy, Government of India, under project no.~12-R\&D-TFR-5.01-0520.

\bibliographystyle{alpha}
\newcommand{\etalchar}[1]{$^{#1}$}


\end{document}